\documentclass[aos,MSNbibl,seceqn,dvips]{arximspdf}
\usepackage{multirow}
\usepackage{dcolumn}
\usepackage{graphicx}

\doi{10.1214/13-AOS1118} 
\volume{41}
\issue{3}
\pubyear{2013}
\firstpage{1350}
\lastpage{1380}

\makeatletter

\let\sv@tabnotetext\tabnotetext
\let\sv@tabnotemark@fmt\tabnotemark@fmt
\long\def\legend#1{{\let\tabnote@indent\leavevmode\sv@tabnotetext[]{}{#1}}}

\newcolumntype{d}[1]{D{.}{.}{#1}}

\newcommand{\rrVert}{\Vert}
\newcommand{\rrvert}{\vert}
\newcommand{\llVert}{\Vert}
\newcommand{\llvert}{\vert}

\newtheorem{lemma}{Lemma}
\newtheorem{theorem}{Theorem}
\newtheorem{proposition}{Proposition}

\newproclaim{remark}{Remark}
\newproclaim{algorithm}{Algorithm}
\newproclaim{assumption}{Assumption}
\newproclaim{model}{Model}

\makeatother

\begin{document}
\begin{frontmatter}

\title{Maximum-likelihood estimation for diffusion processes via
closed-form density expansions}
\runtitle{MLE for diffusion}

\begin{aug}
\author[A]{\fnms{Chenxu} \snm{Li}\corref{}\thanksref{t2}\ead[label=e1]{cxli@gsm.pku.edu.cn}}
\runauthor{C. Li}
\affiliation{Guanghua School of Management, Peking University}
\address[A]{Department of Business Statistics\\
\quad and Econometrics\\
Guanghua School of Management\\
Peking University\\
Beijing, 100871\\
P.R. China\\
\printead{e1}} 
\end{aug}

\thankstext{t2}{This research was supported by the Guanghua School of Management, the
Center for Statistical Sciences, and the Key Laboratory of Mathematical
Economics and Quantitative Finance (Ministry of Education) at Peking
University, as well as the National Natural Science Foundation of China
(Project 11201009).}

\received{\smonth{11} \syear{2011}}
\revised{\smonth{4} \syear{2013}}

%
\begin{abstract}
This paper proposes a widely applicable method of approximate
maxi\-mum-likelihood estimation for multivariate diffusion process from
discretely sampled data. A closed-form asymptotic expansion for
transition density is proposed and accompanied by an algorithm
containing only basic and explicit calculations for delivering any
arbitrary order of the expansion. The likelihood function is thus
approximated explicitly and employed in statistical estimation. The
performance of our method is demonstrated by Monte Carlo simulations
from implementing several examples, which represent a wide range of
commonly used diffusion models. The convergence related to the
expansion and the estimation method are theoretically justified using
the theory of Watanabe [\textit{Ann. Probab.} \textbf{15} (1987) 1--39]
and Yoshida [\textit{J. Japan Statist. Soc.} \textbf{22} (1992)
139--159] on analysis of the generalized random variables under some
standard sufficient conditions.
\end{abstract}

%

\end{frontmatter}

\section{Introduction}
\label{sectionintroduction}

Diffusion processes governed by stochastic differential equations
(hereafter SDE) are widely used in describing the phenomenon of random
fluctuations over time, and even become indispensable for analyzing
high-frequency data; see, for example, Mykland and Zhang
\cite{MyklandZhang2010survey}. Practical application of diffusion
models calls for statistical inference based on discretely monitored
data. The literature has seen a wide spectrum of asymptotically
efficient estimation methods, for example, those based on various
contrast functions proposed in Yoshida~\cite{Yoshida1992Discrete}, %
Kessler~\cite{Kessler1997}, Kessler and S{\o}rensen \cite
{KesslerSorensen1999} and the
references given in S{\o}rensen~\cite{Sorensen2012survey}. Taking the
efficiency, feasibility and generality into account,
maximum-likelihood estimation (hereafter MLE) can be a choice among
others. However, for the increasingly complex real-world dynamics,
likelihood functions (transition densities) are generally not known
in closed-form and thus involve significant challenges in valuation.
This leads to various methods of approximation and the resulting
approximate MLE. The focus of this paper is to propose a widely
applicable closed-form asymptotic expansion for transition density
and thus to apply it in approximate MLE for multivariate diffusion
process.

\subsection{Background}

To approximate likelihood functions, Yoshida \cite
{Yoshida1992Discrete} proposed
to discretize continuous likelihood functions (see, e.g., %
Basawa and Prakasa Rao~\cite{BasawaPrakasaRao1980book}); many others
focused on direct
approximation of likelihood functions (transition densities) for discretely
monitored data, see surveys in, for example, Phillips and Yu \cite
{PhillipsYu2009}, %
Jensen and Poulsen~\cite{JENSENandPOULSEN2002}, Hurn, Jeisman and
Lindsay~\cite{Hurnetal2007} and the references
therein. In particular, among various numerical methods, Lo~\cite{Lo1988MLE}
proposed to employ a numerical solution of Kolmogorov equation for transition
density; Pedersen~\cite{Pedersen1995}, Brandt and Santa-Clara \cite
{BrandtSantaClara2002}, %
Durham and Gallant~\cite{DurhamGallant2002}, Stramer and Yan \cite
{STRAMERandYAN2005}, %
Beskos and Roberts~\cite{BeskosRoberts2005}, %
Beskos et al.~\cite{BeskosPapaspiliopoulosRobertsFearnhead2006}, %
Beskos, Papaspiliopoulos and Roberts~\cite
{BeskosPapaspiliopoulosRoberts2009} and Elerian, Chib and Shephard
\cite{Elerian2001}
advocated the application of various Monte Carlo simulation methods; %
Yu and Phillips~\cite{YuPhillips2001} developed an exact Gaussian
method for
models with a linear drift function; Jensen and Poulsen \cite
{JENSENandPOULSEN2002}
resorted to the techniques of binomial trees. Since all these
numerical methods are computationally demanding, real-world
implementation has necessitated the development of analytical
methods for efficiently approximating transition density. An adhoc
approach is to approximate the model by discretization, for example, the
Euler scheme, and then use the transition density of the discretized
model. Elerian~\cite{Elerian1998} refined such an approximation
via the second order Milstein scheme. Kessler~\cite{Kessler1997} and %
Uchida and Yoshida~\cite{UchidaYoshida2012} employed a more sophisticated
normal-distribution-based approximation via higher order expansions
of the mean and variance.

For approximate MLE of diffusions, %
Dacunha-Castelle and Florens-\break Zmirou~\cite
{DacunhaCastelleandFlorensZmirou1986} is one of the
earliest attempts to apply the idea of small-time expansion of
transition densities, which in principle can be made arbitrarily
accurate. However, their method relies on implicit representation of
moments of Brownian bridge functionals, and thus requires Monte
Carlo simulation in implementation. A milestone is the
ground-breaking work of A{\"{i}}t-Sahalia~\mbox{\cite
{AitSahalia1999JF,AitSahalia2002Econometrica,AitSahalia2008AS}}, which
established the theory of Hermite-polynomial-based analytical
expansion for transition density of diffusion models and the
corresponding approximate MLE. Along the line of
A{\"{i}}t-Sahalia~\cite
{AitSahalia1999JF,AitSahalia2002Econometrica,AitSahalia2008AS}, a
number of substantial refinements and
applications emerged in the literature of likelihood-based
statistical inference (see
surveys in A{\"{i}}t-Sahalia~\cite{AitSahalia2009ARFE}); see, for
example, %
Bakshi and Ju~\cite{BakshiJu2005JB}, Bakshi, Ju and\vadjust{\goodbreak} Ou-Yang \cite
{BakshiJu2006227},
A{\"{\i}}t-Sahalia and Mykland \cite
{AitSahaliaandPerMykland2003Annals,AitSahaliaandPerMykland2003Economstrica},
A{\"{i}}t-Sahalia and Kimmel \cite
{AitSahaliaandKimmel2007SVestimation,AitSahaliaKimmel2002TermStructureJFE},
Li~\cite{Li2010dampeddiffusion},
Egorov, Li and Xu~\cite{EgorovtimeinhomogeneousMLE}, Schaumburg \cite
{SCHAUMBURG2001}, %
A{\"{\i}}t-Sahalia and Yu~\cite{AitSahaliaandYu2006SaddleJE}, Yu
\cite{Yu2007}, %
Filipovi{\'c}, Mayerhofer and Schneider~\cite{Filipovic2011}, Tang and
Chen~\cite{TangChen2009}, %
Xiu~\cite{Xiu2011priceexpansion} and Chang and Chen~\cite{ChangChen2011}.

\subsection{Expansion for likelihood functions and approximate MLE}

Starting from the celebrated Edgeworth expansion for distribution of
standardized summation of independently identically distributed random
variables (see, e.g., Chapter XVI in Feller~\cite{Feller1971}, Chapter
2 in %
Hall~\cite{Hall1995} and Chapter 5 in McCullagh~\cite{MCCULLAGH1987}),
asymptotic expansions have become powerful tools for statistics,
econometrics and many other disciplines in science and technology.
Taking dependence of random variables into account, Mykland
\cite{Mykland1992AoS,Mykland1993,Mykland1994,Mykland1995} established
the theory, calculation and various statistical applications of
martingale expansion, which is further developed in Yoshida
\cite{Yoshida1997martingaleexpansion,Yoshida2001martingaleexpansion}.

Having an analogy with these Edgeworth-type expansions and motivated by MLE
for diffusion processes, I propose a new small-time asymptotic
expansion of
transition density for multivariate diffusions based on the theory of %
Watanabe~\cite{WatanabeAnalysisofWiener1987} and %
Yoshida~\cite{Yoshida1992StatisticsSmallDiffusion,Yoshida1992MLE}.
However, in contrast to the traditional Edgeworth expansions, our
expansion does not require the knowledge of generally implicit
moments, cumulants or characteristic function of the underlying
variable, and thus it is applicable to a wide range of diffusion
processes. Moreover, in analogy
to the verification of validity given in, for example, %
Bhattacharya and Ghosh~\cite{BhattacharyaGhosh1978},
Mykland~\cite{Mykland1992AoS,Mykland1993,Mykland1994,Mykland1995} and
Yoshida \cite
{Yoshida1997martingaleexpansion,Yoshida2001martingaleexpansion} for
Edgeworth type expansions, the
uniform convergence rate (with respect to various parameters) of our density
expansion is proved under some sufficient conditions on the drift and
diffusion coefficients of the underlying diffusion using the theory of %
Watanabe~\cite{WatanabeAnalysisofWiener1987} and %
Yoshida~\cite{Yoshida1992StatisticsSmallDiffusion,Yoshida1992MLE}.
Consequently,
the approximate MLE converges to the true one, and thus inherits its
asymptotic properties. Such results are further demonstrated through
numerical tests and Monte Carlo simulations for some representative examples.

In comparison to the expansion proposed by A{\"{i}}t-Sahalia
\cite{AitSahalia1999JF,AitSahalia2002Econometrica,AitSahalia2008AS},
our method is able to bypass the challenge resulting from the
discussion of reducibility, the
explicity of the \textit{Lamperti} transform (see, e.g., Section 5.2
in %
Karatzas and Shreve~\cite{KaraztasShreve}) and its inversion, as well
as the iterated equations for expressing correction terms, which in
general lead to multidimensional integrals; see Bakshi, Ju and Ou-Yang
\cite{BakshiJu2006227}. Thus it renders an algorithm for practically
obtaining a closed-form expansion (without integrals and implicit
transforms) for transition density up to any arbitrary order, which
serves as a widely applicable tool for approximate MLE. Even after the
\textit{Lamperti} transform, our expansion employs a completely
different nature comparing with those proposed in A{\"{i}}t-Sahalia
\cite{AitSahalia1999JF,AitSahalia2002Econometrica,AitSahalia2008AS},
which hinge on expansions in an orthogonal basis consisting of Hermite
polynomials and expansions of each coefficient expressed by an
expectation of a smooth functional of the transformed
variable via an iterated Dynkin formula; see Section~4 in %
A{\"{\i}}t-Sahalia~\cite{AitSahalia2002Econometrica}.\vadjust{\goodbreak}

Moreover, our method is different from the existing theory of
large-devia\-tions-based expansions, which were discussed in, for example,
Azencott~\cite{Azencott1984}, %
Bismut~\cite{Bismut1984}, Ben Arous~\cite{BenArous1988} and L{\'e}andre
\cite{Leandre1992}, and
given probabilistic representation in %
Watanabe~\cite{WatanabeAnalysisofWiener1987} for the purpose of
investigating the analytical structure of heat kernel in
differential geometry. Large-deviations-based asymptotic expansions
involve Riemannian distance (implied by the true but generally
unknown transition density) and higher order correction terms.
Except for some special cases, they rarely admit closed-form
expressions by solving the corresponding variational problems.
However, for practical implementation of statistical estimation,
relatively simple closed-form approximations are usually favorable.

The rest of this paper is organized as follows. In Section \ref
{SectionModelMLE}, the model is introduced with some technical assumptions
and the maximum-likelihood estimation problem is formulated. In Section
\ref{sectionexpansionframework}, the transition density expansion is proposed
with closed-form correction terms of any arbitrary order for general
multivariate diffusion processes, and the uniform convergence of
the expansion is \mbox{established}. In Section \ref
{sectionimplementationexamples}, numerical
performance of the density expansion is demonstrated through examples. In
Section~\ref{SectionAMLE}, the asymptotic properties of the consequent
approximate MLE are established. In Section~\ref{SectionAMLEMC}, Monte
Carlo evidence for the approximate MLE is provided. In Section \ref
{sectionconcludingremarks}, the paper is concluded and some
\mbox{opportunities} for future research are outlined. Appendix
\ref{subsectionCE} provides an algorithm for explicitly calculating a
type of conditional expectation, which plays an important role in the
closed-form expansion. Appendix~\ref{appendixproofs} contains all
proofs. The supplementary material~\cite{autokey73} collects some
concrete formulas for illustration, figures for exhibiting detailed
numerical performance, additional and alternative output of simulation
results, more examples, a brief introduction to the theory of
Watanabe--Yoshida and the proof of a technical lemma.

\section{The model and maximum-likelihood estimation}
\label{SectionModelMLE}

Assuming known parametric form of the drift vector function
$\mu=(\mu_{1},\ldots,\mu_{m})\dvtx \mathbb{R}^{m}\rightarrow \mathbb{R}^{m}$
and the dispersion matrix $\sigma=(\sigma_{ij})_{m\times d}$:
$%
\mathbb{R}^{m}\rightarrow\mathbb{R}^{m\times d}$ with unknown
parameter $%
\theta$ belonging to a compact set $\Theta\subset\mathbb{R}^{k}$,
an $m$-dimensional time-homogenous diffusion $X$ is modeled by an SDE,
%
\begin{equation}\label{GeneralmodelX}
dX(t)=\mu\bigl(X(t);\theta\bigr)\,dt+\sigma\bigl(X(t);\theta\bigr)\,dW(t),\qquad
X(0)=x_{0},
\end{equation}
where $\{W(t)\}$ is a $d$-dimensional standard Brownian motion. Let $%
E\subset\mathbb{R}^{m}$ denote the state space of $X$. Without loss
of generality, we assume $m=d$ throughout the paper.

By the time-homogeneity nature of diffusion $X$, let $p_{X}(\Delta,x|x_{0};\theta)$ denote the conditional density of $X(t+\Delta)$
given $%
X(t)=x_{0}$, that is,
\[
\mathbb{P}\bigl(X(t+\Delta)\in dx|X(t)=x_{0}\bigr)=p_{X}(
\Delta,x|x_{0};\theta)\,dx.
\]
Based on the discrete observations of $X$ at time grids
$\{\Delta,2\Delta,\ldots,n\Delta\}$, which correspond to the daily,
weekly or monthly monitoring, etc., the likelihood function is
constructed as
%
\begin{equation}\label{likelihood}
l_{n}(\theta)=\prod_{i=1}^{n}p_{X}
\bigl(\Delta t,X(i\Delta)|X\bigl((i-1)\Delta\bigr);\theta\bigr);
\end{equation}
the corresponding log-likelihood function admits the following form:
%
\begin{equation}\label{loglikelihood}
\ell_{n}(\theta)=\sum_{i=1}^{n}L_{i}(
\theta),
\end{equation}
where the $\log$ transition density is
%
\begin{equation}\label{logdensity}
L_{i}(\theta)=\log\bigl[ p_{X}\bigl(\Delta,X(i\Delta)|X
\bigl((i-1)\Delta\bigr);\theta\bigr)%
\bigr].
\end{equation}
Maximum-likelihood estimation is to identify the optimizer in $\theta
\in
\Theta$ for (\ref{likelihood}) or equivalently (\ref{loglikelihood}).
However, except for some simple models, (\ref{likelihood}) and (\ref
{loglikelihood}) rarely admit closed-form expressions.

For ease of exposition, we introduce some technical assumptions. Let $%
A(x;\theta)=\sigma(x;\theta)\sigma(x;\theta)^{T}$ denote the diffusion
matrix.

\begin{assumption}
\label{assumptionpositivedefinite} The diffusion matrix
$A(x;\theta)$ is positive definite, that is, $\det A(x;\theta)>0$, for
any $(x,\theta)\in E\times\Theta$.
\end{assumption}

\begin{assumption}
\label{aussumptionboundedderivatives} For each integer $k\geq1$, the
$k$th order derivatives in $x$ of the functions $\mu(x;\theta)$ and
$\sigma(x;\theta)$ exist, and they are uniformly bounded for any
$(x,\theta)\in E\times\Theta$.
\end{assumption}

\begin{assumption}
\label{assumption3timesdifferentiable} The transition density
$p_{X}(\Delta,x|x_{0};\theta)$ is continuous in $\theta\in\Theta$,
and the log-likelihood function (\ref{loglikelihood}) admits a
unique maximizer in the parameter set $\Theta$.
\end{assumption}

Assumptions~\ref{assumptionpositivedefinite} and \ref
{aussumptionboundedderivatives} are conventionally proposed in the
study of stochastic differential equations; see, for example, Ikeda and
Watanabe~\cite{IkedaWatanabe1989}. They are sufficient (but not
necessary) to guarantee the existence and uniqueness of the solution
and other desirable technical properties. For convenience, the
theoretical proofs given in Appendix~\ref{appendixproofs} are based on
these conditions. However, as is shown in Sections \ref
{sectionimplementationexamples} and~\ref{SectionAMLEMC}, numerical
examples suggest that the method proposed in this paper is applicable
to a wide range of commonly used models, rather than confined to those
strictly satisfying these sufficient (but not necessary) conditions.
Assumption~\ref{assumption3timesdifferentiable} collects two standard
conditions for maximum likelihood estimation. In particular, for the
continuity (and higher differentiability) of the transition density in
the parameter, sufficient conditions based on the smoothness of the
drift and dispersion functions can be found in, for example, Azencott
\cite{Azencott1984} and A{\"{\i}}t-Sahalia \cite
{AitSahalia2002Econometrica}. Theoretical relaxation of these
conditions may involve case-by-case treatment and standard
approximation argument, which is beyond the scope of this paper and can
be regarded as a future research topic.

\section{A closed-form expansion for transition density}
\label{sectionexpansionframework}

The method of approximate maximum-likelihood estimation proposed in
this paper relies on a closed-form expansion for transition density of
any arbitrary diffusion process. \mbox{Bypassing} the discussion of the
\textit{Lamperti} transform and the reducibility issue as in
A{\"{i}}t-Sahalia
\cite{AitSahalia1999JF,AitSahalia2002Econometrica,AitSahalia2008AS},
our starting point stands on the fact that the transition density can
be expressed as
%
\begin{equation}\label{ExpectationDiracDensity}
p_{X}(\Delta,x|x_{0};\theta)=\mathbb{E} \bigl[
\delta\bigl(X(\Delta)-x\bigr)|X(0)=x_{0} \bigr],
\end{equation}
where $\delta(z)$ is the Dirac Delta function centered at $0$ for
some variable $z$. More precisely, $\delta(z)$ is defined as a
generalized function (distribution) such that it is zero for all
values of $z$ except when it is zero, and its integral from $-\infty
$ to $\infty$ is equal to one; see, for example,
Kanwal~\cite{Kanwal2004generalizedfunctions} for more details.
Watanabe~\cite{WatanabeAnalysisofWiener1987} established the validity of
(\ref{ExpectationDiracDensity}) through the theory of generalized
random variables and expressed correction terms of
large-deviations-based density expansion as implicit expectation
forms by separately treating the cases of diagonal ($x=x_{0}$) and
off-diagonal ($x\neq x_{0}$). In particular, the off-diagonal
($x\neq x_{0}$) expansion depends on a generally implicit
variational formulation for Riemanian distance. From the viewpoint
of statistical applications where $X(\Delta)\neq X(0)$
(corresponding to $x\neq x_{0}$) happens almost surely, the
expansion proposed in Watanabe~\cite{WatanabeAnalysisofWiener1987} is
impractical due to
high computational costs. In the literature of statistical inference, %
(\ref{ExpectationDiracDensity}) has been employed in
Pedersen~\cite{Pedersen1995} for simulation-based approximate MLE. In this
section, we propose a new expansion of the transition density which
universally treats the diagonal ($x=x_{0}$) and off-diagonal ($x\neq
x_{0}$) cases. Heuristically speaking, our method hinges on a
Taylor-like expansion of a standardized version of $\delta(X(\Delta
)-x)$, which results in closed-form formulas for any arbitrary
correction term.

\subsection{Basic setup and notation}
\label{subsectionassumptionbasicsetup}

Let $\epsilon=\sqrt{\Delta}$ be a small parameter based on which
an asymptotic expansion is carried out. By rescaling the model (\ref
{GeneralmodelX}) to bring forth finer local behavior of the
diffusion process, we let $X^{\epsilon}(t):=X(\epsilon^{2}t)$.
Integral substitution and the Brownian scaling property yield that
%
\begin{equation}\label{scaleddiffusion}\quad
dX^{\epsilon}(t)=\epsilon^{2}\mu\bigl(X^{\epsilon}(t);\theta
\bigr)\,dt+\epsilon\sigma\bigl(X^{\epsilon}(t);\theta\bigr)\,dW^{\epsilon
}(t),\qquad X^{\epsilon}(0)=x_{0},
\end{equation}
where $\{W^{\epsilon}(t)\}$ is a $m$-dimensional standard Brownian
motion. For notation simplicity, we write the scaled Brownian motion
$W^{\epsilon}(t) $ as $W(t)$ and drop the parameter $\theta$ in
what follows.

Let us introduce a vector function
$b(x)=(b_{1}(x),b_{2}(x),\ldots,b_{m}(x))^{T}$ defined by
%
\begin{equation}\label{b}
b_{i}(x)=\mu_{i}(x)-\frac{1}{2}\sum
_{k=1}^{m}\sum_{j=1}^{m}
\sigma_{kj}(x)\,
\frac{\partial}{\partial x_{k}}\sigma_{ij}(x)
\end{equation}
and construct the following differential operators:
%
\begin{equation}\label{A0AjOperators}\qquad
\mathcal{A}_{0}:=\sum_{i=1}^{m}b_{i}(x)\,
\frac{\partial}{\partial
x_{i}}\quad\mbox{and}\quad\mathcal{A}_{j}:=\sum
_{i=1}^{m}\sigma_{ij}(x)\,\frac{\partial}{\partial
x_{i}}\qquad
\mbox{for }j=1,\ldots,m,
\end{equation}
which map vector-valued functions to vector-valued functions of the
same dimension, respectively. More precisely, for any
$\nu\in\mathbb{N}$ and a $\nu$-dimensional vector-valued function
$\varphi(x)=(\varphi
_{1}(x),\varphi_{2}(x),\ldots,\varphi_{\nu}(x))^{T}$,
\[
\bigl( \mathcal{A}_{0} ( \varphi) \bigr) (x)= \Biggl( \sum
_{i=1}^{m}b_{i}(x)\,\frac{\partial\varphi_{1}(x)}{\partial
x_{i}},\sum_{i=1}^{m}b_{i}(x) \,\frac{\partial\varphi_{2}(x)}{\partial
x_{i}},\ldots,\sum_{i=1}^{m}b_{i}(x) \,\frac{\partial\varphi_{\nu
}(x)}{\partial x_{i}} \Biggr) ^{T}
\]
and
\[
\bigl( \mathcal{A}_{j} ( \varphi) \bigr) (x)= \Biggl( \sum
_{i=1}^{m}\sigma_{ij}(x)\,
\frac{\partial\varphi_{1}(x)}{\partial
x_{i}},\sum_{i=1}^{m}
\sigma_{ij}(x)\,\frac{\partial\varphi
_{2}(x)}{\partial x_{i}},\ldots,\sum
_{i=1}^{m}\sigma_{ij}(x)\,\frac{\partial\varphi_{\nu
}(x)}{%
\partial x_{i}}
\Biggr) ^{T}
\]
for $j=1,2,\ldots,m$.

For an index $\mathbf{i}=(i_{1},\ldots,i_{n})\in\{0,1,2,\ldots,m\}^{n}$
and a right-continuous sto\-chastic process $\{f(t)\}$, define an iterated
Stratonovich integral with integrand $f$ as
%
\begin{equation}\label{StratIntfgeneral}\qquad
J_{\mathbf{i}}[f](t):=\int_{0}^{t}\int
_{0}^{t_{1}}\cdots\int_{0}^{t_{n-1}}f(t_{n})
\circ dW_{i_{n}}(t_{n})\cdots\circ dW_{i_{2}}(t_{2})
\circ dW_{i_{1}}(t_{1}),
\end{equation}
where $\circ$ denotes stochastic integral in the Stratonovich sense. Note
that $J_{\mathbf{i}}[f](t)$ is recursively defined from inside to
outside; see page 174 of
Kloeden and Platen~\cite{KloedenPlaten1999}. For ease of exposition, the
order of iterated integrations defined in this paper is the reverse of that
in Kloeden and Platen~\cite{KloedenPlaten1999} for any arbitrary
index. To lighten the
notation, for $f\equiv1$, the integral $J_{\mathbf{i}}[1](t)$ is
abbreviated to $J_{\mathbf{i}}(t)$. By convention, let $W_{0}(t):=t$ and
define
%
\begin{equation}\label{inorm}
\llVert\mathbf{i}\rrVert:=\sum_{k=1}^{n}
[ 2\cdot1_{\{i_{k}=0\}}+1_{\{i_{k}\neq0\}} ]
\end{equation}
as a ``norm'' of index $\mathbf{i}$, which
counts an index $k$ with $i_{k}=0$ twice.

By viewing $X^{\epsilon}(1)$ as a function of $\epsilon$, it is natural
to obtain a pathwise expansion in $\epsilon$ with random coefficients,
which serves as a foundation for our transition density expansion. According
to Watanabe~\cite{WatanabeAnalysisofWiener1987}, I introduce the following
coefficient function $C_{\mathbf{i}}(x_{0})$ defined by iterative
application of the differential operators (\ref{A0AjOperators}):
%
\begin{equation}\label{Ccoefficient}
C_{\mathbf{i}}(x_{0}):=\mathcal{A}_{i_{n}} \bigl( \cdots
\bigl( \mathcal{A}%
_{i_{3}} \bigl( \mathcal{A}_{i_{2}} (
\sigma_{\cdot
i_{1}} ) \bigr) \bigr) \cdots\bigr) (x_{0})
\end{equation}
for an index $\mathbf{i}=(i_{1},\ldots,i_{n})$. Here, for $i_{1}\in
\{1,2,\ldots,m\}$, the vector $\sigma_{\cdot i_{1}}(x)=(\sigma
_{1i_{1}}(x),\ldots,\sigma_{mi_{1}}(x))^{T}$ denotes the $i_{1}$th column
vector of the dispersion matrix $\sigma(x)$, for $i_{1}=0$, $\sigma
_{\cdot
0}(x)$ refers to the vector $b(x)$ defined in (\ref{b}).

Using vector function (\ref{b}), the scaled diffusion (\ref
{scaleddiffusion}) can be equivalently expressed as the following
stochastic differential equation in the Stratonovich sense (see,
e.g., Section 3.3 in Karatzas and Shreve~\cite{KaraztasShreve}), that is,
\[
dX^{\epsilon}(t)=\epsilon^{2}b\bigl(X^{\epsilon}(t)\bigr)\,dt+
\epsilon\sigma\bigl(X^{\epsilon}(t)\bigr)\circ dW(t).
\]
Thus, similarly to Theorem 3.3 in Watanabe \cite
{WatanabeAnalysisofWiener1987},
it is easy to obtain a closed-form pathwise expansion of $X^{\epsilon}(1)$
from successive applications of the It\^{o} formula.

\begin{lemma}
$X^{\epsilon}(1)$ admits the following pathwise asymptotic expansion:
%
\begin{equation}\label{Xpathwiseexpansion}
X^{\epsilon}(1)=\sum_{k=0}^{J}F_{k}
\epsilon^{k}+\mathcal{O}\bigl(\epsilon^{J+1}\bigr)
\end{equation}
for any $J\in N$. Here, $F_{0}=x_{0}$ and $F_{k}$ can be written as a
closed-form linear combination of iterated Stratonovich integrals, that
is,
%
\begin{equation}\label{Fk}
F_{k}=\sum_{\llVert\mathbf{i}\rrVert=k}C_{\mathbf
{i}}(x_{0})J_{%
\mathbf{i}}(1)
\end{equation}
for $k=1,2,\ldots\,$, where the integral $J_{\mathbf{i}}(1)$, the norm $
\llVert\mathbf{i}\rrVert$ and coefficient $C_{\mathbf{i}}(x_{0})$
are defined in (\ref{StratIntfgeneral}), (\ref{inorm}) and (\ref
{Ccoefficient}), respectively.
\end{lemma}

For any arbitrary dimension $r=1,2,\ldots,m$, one has the element-wise
form of the expansion (\ref{Xpathwiseexpansion}) as
$X_{r}^{\epsilon}(1)=\sum_{k=0}^{J}F_{k,r}\epsilon
^{k}+\mathcal{O}(\epsilon^{J+1})$ where
%
\begin{eqnarray}\label{Ccoefficientsr}
F_{k,r}&=&\sum_{\llVert\mathbf{i}\rrVert
=k}C_{\mathbf{i},r}(x_{0})J_{%
\mathbf{i}}(1)
\end{eqnarray}
with
\[
C_{\mathbf{i},r}(x_{0}):=\mathcal{A}_{i_{n}}
\bigl( \cdots\bigl( \mathcal{A}%
_{i_{3}} \bigl(
\mathcal{A}_{i_{2}} ( \sigma_{ri_{1}} ) \bigr) \bigr)\cdots\bigr)
(x_{0})
\]
for $\mathbf{i}=(i_{1},\ldots,i_{n})$. Note that (\ref
{Xpathwiseexpansion}%
) is different from the Wiener chaos decomposition (see, e.g., %
Nualart~\cite{NualartMCbook}), which employs an alternative way of representing
random variables. The validity of the pathwise expansion (\ref
{Xpathwiseexpansion}) and other expansions introduced in the next
subsection can be
rigorously guaranteed by the theory of %
Watanabe~\cite{WatanabeAnalysisofWiener1987} and %
Yoshida~\cite{Yoshida1992StatisticsSmallDiffusion,Yoshida1992MLE}. For
ease of
exposition, we focus on the derivation of density expansion in this and the
following subsection and articulate the validity issue in Section
\ref{Subsectionconvergencedensityexpansion}.\vadjust{\goodbreak}

We introduce an $m$-dimensional correlated Brownian motion
%
\begin{equation}\label{BMB}
B(t)=\bigl(B_{1}(t),B_{2}(t),\ldots,B_{m}(t)
\bigr) \qquad\mbox{with } B_{k}(t)=\frac
{\sum_{j=1}^{m}\sigma_{kj}(x_{0})W_{j}(t)}{\sqrt{%
\sum_{j=1}^{m}\sigma_{kj}^{2}(x_{0})}}\hspace*{-32pt}
\end{equation}
for $k=1,2,\ldots,m$. Thus, the leading term $F_{1}$ can be expressed
as
\[
F_{1}= \Biggl( \sqrt{\sum_{j=1}^{m}
\sigma_{1j}^{2}(x_{0})}%
B_{1}(1),\sqrt{\sum_{j=1}^{m}
\sigma_{2j}^{2}(x_{0})}%
B_{2}(1),\ldots,\sqrt{\sum_{j=1}^{m}
\sigma_{mj}^{2}(x_{0})}%
B_{m}(1) \Biggr).
\]
Let $D(x)$ be a diagonal matrix defined by
%
\begin{equation}\label{Dx}
D(x):=\operatorname{diag} \biggl(
\frac{1}{\sqrt{\sum_{j=1}^{m}\sigma_{1j}^{2}(x)}},\frac{1}{\sqrt{\sum_{j=1}^{m}
\sigma_{2j}^{2}(x)}},\ldots,\frac{1}{\sqrt{\sum_{j=1}^{m}\sigma_{mj}^{2}(x)}}
\biggr).\hspace*{-32pt}
\end{equation}
It follows that $B(t)=D(x_{0})\sigma(x_{0})W(t)$ and
$D(x_{0})F_{1}=B(1)$. Furthermore, the correlation of $B_{k}(t)$ and
$B_{l}(t)$ for $k\neq l$ is given by
\[
\rho_{kl}(x_{0}):=\operatorname{Corr}\bigl(B_{k}(t),B_{l}(t)
\bigr)=\frac{\sum_{j=1}^{m}\sigma_{kj}(x_{0})\sigma_{lj}(x_{0})}{\sqrt{\sum_{j=1}^{m}%
\sigma_{kj}^{2}(x_{0})}\sqrt{\sum_{j=1}^{m}
\sigma_{lj}^{2}(x_{0})}}.
\]
So, the covariance matrix of $B(1)$ is
%
\begin{equation}\label{CovariancematrixB}
\Sigma(x_{0})=\bigl(\rho_{ij}(x_{0})
\bigr)_{m\times m}=D(x_{0})\sigma(x_{0})\sigma
(x_{0})^{T}D(x_{0}).
\end{equation}
It follows that Assumption~\ref{assumptionpositivedefinite} is equivalent
to the positive definite property of the correlation matrix $\Sigma(x_{0})$
and the nonsingularity of the dispersion matrix $\sigma(x_{0})$, that
is, $%
\det A(x_{0})>0\Longleftrightarrow\det\Sigma
(x_{0})>0\Longleftrightarrow
\det\sigma(x_{0})>0$. Finally, for any index
$i\in
\{1,2,\ldots,m\}$ and differentiable function $u(y)$ with $y\in
\mathbb{R}%
^{m}$, we introduce the following differential operator:
%
\begin{equation}\label{Doperators}
\mathcal{D}_{i}u(y):=\frac{\partial u(y)}{\partial y_{i}}-u(y) \bigl
(\Sigma
(x_{0})^{-1}y\bigr)_{i},
\end{equation}
where $(\Sigma(x_{0})^{-1}y)_{i}$ denotes the $i$th element of the
vector $\Sigma(x_{0})^{-1}y$.

\subsection{Asymptotic expansion for transition densities: A general
framework}

Employing the scaled diffusion $X^{\epsilon}(t)=X(\epsilon^{2}t)$
with $%
\epsilon=\sqrt{\Delta}$, the expectation representation (\ref
{ExpectationDiracDensity}) for transition density can be expressed
as
%
\begin{equation}\label{pEdelta}
p_{X}(\Delta,x|x_{0};\theta)=\mathbb{E} \bigl[
\delta\bigl(X^{\epsilon
}(1)-x\bigr)|X^{\epsilon}(0)=x_{0}
\bigr].
\end{equation}
To guarantee the convergence, our expansion procedure begins with
standardizing $X^{\epsilon}(1)$ to
%
\begin{equation}\label{standardizationXY}
Y^{\epsilon}:=D(x_{0})\frac{X^{\epsilon}(1)-x_{0}}{\epsilon
}=D(x_{0})%
\frac{X^{\epsilon}(1)-x_{0}}{\sqrt{\Delta}},
\end{equation}
which converges to a nonconstant random variable (a multivariate
normal in
our case), see Watanabe~\cite{WatanabeAnalysisofWiener1987} and %
Yoshida~\cite{Yoshida1992StatisticsSmallDiffusion,Yoshida1992MLE} for a
similar setting. Indeed, based on the Brownian motion defined in
(\ref{BMB}) and the fact $D(x_{0})F_{1}=B(1)$, the $j$th component
of $Y^{\epsilon}(1)$ satisfies that
%
\begin{equation}\label{Y}
Y_{j}^{\epsilon}:=\frac{X_{j}^{\epsilon}(1)-x_{0j}}{\epsilon\sqrt{%
\sum_{i=1}^{d}\sigma_{ji}^{2}(x_{0})}}\rightarrow
B_{j}(1)\qquad\mbox{as }\epsilon\rightarrow0
\end{equation}
for $j=1,2,\ldots,m$. It is worth noting that
Watanabe~\cite{WatanabeAnalysisofWiener1987} employed an alternative
standardization method (see Theorem 3.7 in %
Watanabe~\cite{WatanabeAnalysisofWiener1987}) in constructing the
implicit expectation representation for the correction terms of
large-deviations-based
density expansion for the case of $x\neq x_{0}$; see Theorem 3.8 in %
Watanabe~\cite{WatanabeAnalysisofWiener1987}.

Owing to (\ref{standardizationXY}), the pathwise expansion (\ref
{Xpathwiseexpansion}) implies that
%
\begin{equation}\label{Yexpansion}
Y^{\epsilon}=\sum_{i=0}^{J}Y_{i}
\epsilon^{i}+\mathcal{O}\bigl(\epsilon^{J+1}\bigr)
\qquad\mbox{with
}Y_{i}=D(x_{0})F_{i+1}
\end{equation}
for any $J\in\mathbb{N}$. Thus, based on (\ref{pEdelta}%
), a Jacobian transform resulting from the change of variable in (\ref
{standardizationXY}) yields the following representation of the
density of $X^{\epsilon}(1)$ based on that of $Y^{\epsilon}$, that is,
\[
p_{X}(\Delta,x|x_{0};\theta)= \biggl( \frac{1}{\sqrt{\Delta}}
\biggr) ^{m}\det D(x_{0})\mathbb{E} \bigl[ \delta
\bigl(Y^{\epsilon
}-y\bigr)|X(0)=x_{0} \bigr],
\]
where $y=D(x_{0}) ( x-x_{0} ) /\sqrt{\Delta}$. For ease of
exposition, the initial condition $X(0)=x_{0}$ is omitted in what
follows. So, the key task is to develop an asymptotic expansion for $%
\mathbb{E} [ \delta(Y^{\epsilon}-y) ] $ around $\epsilon=0$.

Based on the theory of Watanabe~\cite{WatanabeAnalysisofWiener1987}
and %
Yoshida~\cite{Yoshida1992StatisticsSmallDiffusion,Yoshida1992MLE}, the
Dirac Delta function can be manipulated as a function for many
purposes, though it can be formally defined as a distribution. Based
on the expansion of $Y^{\epsilon}$ and heuristic application of
classical rule for differentiating composite functions [the Dirac
Delta function $\delta(\cdot-y)$ acting on $Y^{\epsilon}$ as a
function of $\epsilon$], one is able to obtain a Taylor expansion
of $\delta(Y^{\epsilon}-y)$ as
%
\begin{equation}\label{DeltaYexpansion}
\delta\bigl(Y^{\epsilon}-y\bigr)=\sum_{k=0}^{J}
\Phi_{k}(y)\epsilon^{k}+\mathcal{O}%
\bigl(\epsilon^{J+1}\bigr)
\end{equation}
for any $J\in\mathbb{N}$, where $\Phi_{k}(y)$ represents the
coefficient of the $k$th expansion term. Thus, the following expansion
is immediately implied:
%
\begin{equation}\label{EdeltaY}
\mathbb{E} \bigl[ \delta\bigl(Y^{\epsilon}-y\bigr) \bigr]:=\sum
_{k=0}^{J}\Omega_{k}(y)
\epsilon^{k}+\mathcal{O}\bigl(\epsilon^{J+1}\bigr),
\end{equation}
where $\Omega_{k}(y):=\mathbb{E}\Phi_{k}(y)$ will be explicitly
derived and the remainder term is interpreted in the sense of
classical calculus. Thus, the approximate transition density for $X$
up to the $J$th order is proposed as
%
\begin{eqnarray}\label{PXJ}\quad
p_{X}^{(J)}(\Delta,x|x_{0};\theta)&:=& \biggl(
\frac{1}{\epsilon
} \biggr) ^{m}\det D(x_{0})\sum
_{k=0}^{J}\Omega_{k} \biggl(
D(x_{0})\frac
{x-x_{0}}{%
\epsilon} \biggr) \epsilon^{k}
\nonumber\\[-8pt]\\[-8pt]
&=& \biggl( \frac{1}{\sqrt{\Delta}} \biggr) ^{m}\det D(x_{0})
\sum_{k=0}^{J}\Omega_{k} \biggl(
D(x_{0})\frac{x-x_{0}}{\sqrt{\Delta}%
} \biggr) \Delta^{{k}/{2}}.\nonumber
\end{eqnarray}
The convergence of this expansion (guaranteed by the theory of %
Watanabe~\cite{WatanabeAnalysisofWiener1987} and %
Yoshida~\cite{Yoshida1992StatisticsSmallDiffusion,Yoshida1992MLE})
will be
discussed in Section~\ref{Subsectionconvergencedensityexpansion}.

As outlined in the whole framework, our idea naturally originates from
pathwise expansion of a standardized random variable. However, explicit
calculation of the correction terms $\Omega_{k}$ is still a challenging
issue. In what follows, we will give a general closed-form formula.
Based on (\ref{Y}), (\ref{Yexpansion}), (\ref{DeltaYexpansion}) and
(\ref{EdeltaY}), it is straightforward to
find the
leading term as
%
\begin{eqnarray}\label{leadingordergeneral}
\Omega_{0}(y)&=&\mathbb{E} \bigl[ \delta(Y_{0}-y) \bigr] =
\mathbb{E} \bigl[ \delta\bigl(B(1)-y\bigr) \bigr] \nonumber\\[-8pt]\\[-8pt]
&=&\phi_{\Sigma(x_{0})}(y):=
\frac{\exp(
-y^{T}\Sigma(x_{0})^{-1}y/2 ) }{(2\pi)^{{m}/{2}}(\det
\Sigma
(x_{0}))^{{1}/{2}}},\nonumber
\end{eqnarray}
where $\Sigma(x_{0})$ is defined in (\ref{CovariancematrixB}).

To express $\Omega_{k}(y)$ for arbitrary $k\in\mathbb{N}$, we
introduce an
index set
%
\begin{eqnarray}\label{Sset}
S_{k}&=&\bigl\{ \bigl( l,\mathbf{r}(l),\mathbf{j}(l) \bigr)
|l=1,2,\ldots, \mathbf{r}(l)=(r_{1},r_{2},\ldots,r_{l})\in\{1,2,\ldots,m\}^{l},\hspace*{-27pt}
\nonumber\\[-8pt]\\[-8pt]
&&\hspace*{24.5pt}\mathbf{j}(l)=(j_{1},j_{2},\ldots,j_{l})
\mbox{ with }j_{i}\geq1\mbox{ and
}j_{1}+j_{2}+\cdots+j_{l}=k\bigr\}.\hspace*{-27pt}\nonumber
\end{eqnarray}
As building blocks, let
$P_{(\mathbf{i}_{1},\mathbf{i}_{2},\ldots,\mathbf{i}_{l})}(z)$ denote a
multivariate function in $z=(z_{1},z_{2},\ldots,z_{m})\in\mathbb{R}^{m}$ defined by the conditional expectation of
multiplication of iterated Stratonovich integrals with arbitrary
indices $%
\mathbf{i}_{1},\mathbf{i}_{2},\ldots,\mathbf{i}_{l}$, that is,
%
\begin{equation}\label{ce_product_stratonovichintegrals}
P_{(\mathbf{i}_{1},\mathbf{i}_{2},\ldots,\mathbf
{i}_{l})}(z):=\mathbb{E}%
\Biggl( \prod
_{\omega=1}^{l}J_{\mathbf{i}_{\omega
}}(1)|W(1)=z \Biggr),
\end{equation}
which can be explicitly calculated as a multivariate polynomial
according to an effective algorithm proposed in Appendix
\ref{subsectionCE}.

Now, we will give an explicit formula for obtaining any arbitrary
correction term $%
\Omega_{k}(y)$ under any arbitrary multivariate diffusion process in the
following proposition, which can be implemented using only basic and
explicit calculations in any symbolic software package, for example,
Mathematica.

\begin{theorem}
\label{propositiongeneralOmegak} For any $k\in\mathbb{N}$, the correction
term $\Omega_{k}(y)$ in (\ref{PXJ}) admits the following explicit
expression:
%
\begin{equation}\label{formulageneralOmegak}
\Omega_{k}(y)= \biggl( \sum_{ ( l,\mathbf{r}(l),\mathbf
{j}(l) ) \in
S_{k}}Q_{ ( l,\mathbf{r}(l),\mathbf{j}(l) ) }(y)
\biggr) \phi_{\Sigma(x_{0})}(y),
\end{equation}
where $Q_{ ( l,\mathbf{r}(l),\mathbf{j}(l)) }(y)$
is a polynomial explicitly calculated from
%
\begin{eqnarray}\label{Qlrj}
&&Q_{ ( l,\mathbf{r}(l),\mathbf{j}(l))
}(y) \nonumber\\
&&\qquad=%
\frac{(-1)^{l}}{l!}\sum
_{ \{ (\mathbf{i}_{1},\mathbf
{i}_{2},\ldots,%
\mathbf{i}_{l})|\llVert\mathbf{i}_{\omega}\rrVert
=j_{\omega
}+1,\omega=1,2,\ldots,l \} }\prod_{\omega=1}^{l}
\bigl[ C_{%
\mathbf{i}_{\omega},r_{\omega}}(x_{0})D_{r_{\omega}r_{\omega
}}(x_{0})%
\bigr]
\nonumber
\\
&&\qquad\quad{}\times\mathcal{D}_{r_{1}} \bigl( \mathcal{D}_{r_{2}} \bigl(\cdots\mathcal{D%
}_{r_{l}} \bigl( P_{(\mathbf{i}_{1},\mathbf{i}_{2},\ldots,\mathbf{i}
_{l})} \bigl(
\sigma(x_{0})^{-1}D(x_{0})^{-1}y \bigr)
\bigr) \cdots\bigr) \bigr)
\end{eqnarray}
for the index $ ( l,\mathbf{r}(l),\mathbf{j}(l) ) = (
l,(r_{1},r_{2},\ldots,r_{l}),(j_{1},j_{2},\ldots,j_{l}) ) \in S_{k}$.
Here, $S_{k}$, $\phi_{\Sigma(x_{0})}(y)$, $\llVert\cdot
\rrVert$, $C_{\mathbf{i}_{\omega},r_{\omega}}(x_{0})$,
$D_{r_{\omega
}r_{\omega}}(x_{0})$, $\mathcal{D}_{r_{i}}$ and $P_{(\mathbf{i}_{1},%
\mathbf{i}_{2},\ldots,\mathbf{i}_{l})} ( \cdot) $ are
defined in
(\ref{Sset}), (\ref{leadingordergeneral}), (\ref{inorm}), (\ref
{Ccoefficient}), (\ref{Dx}), (\ref{Doperators}) and (\ref
{ce_product_stratonovichintegrals}), respectively.
\end{theorem}

\begin{pf}
See Appendix~\ref{appendixproofs}.
\end{pf}

An algorithm for explicitly calculating conditional expectation
(\ref{ce_product_stratonovichintegrals}), which plays an important
role in completing the closed-form correction terms as proposed in
Theorem~\ref{propositiongeneralOmegak}, is given in %
Appendix~\ref{subsectionCE}. Regardless of the dimension of
diffusion
processes, I concretely exemplify the closed-form expression (\ref
{formulageneralOmegak}) by the first three correction terms in
the supplementary material~\cite{autokey73}. With $\Omega_{k}$ given by (\ref
{formulageneralOmegak}), a closed-form expansion for transition density
can be constructed via (\ref{PXJ}). 

\subsection{Convergence of the expansion}
\label{Subsectionconvergencedensityexpansion}

In this subsection, we establish the uniform convergence of the
asymptotic expansion (\ref{PXJ}), which will serve as an important
building block for the asymptotic properties of approximate
maximum-likelihood estimation discussed in Section~\ref{SectionAMLE}.
Theoretically speaking, unlike the Hermite-polynomial-based method in
A{\"{i}}t-Sahalia
\cite{AitSahalia1999JF,AitSahalia2002Econometrica,AitSahalia2008AS},
which allows justification of convergence as more correction terms are
added, our new method is a Taylor-like asymptotic expansion, which is
established in the neighborhood of $\Delta=0$. However, as demonstrated
in the numerical experiments and Monte Carlo evidence in Sections
\ref{sectionimplementationexamples} and \ref {SectionAMLEMC},
respectively, accuracy of the expansion is enhanced as $J$ increases
while holding $\Delta$ fixed. Based on the theory of
Watanabe~\cite{WatanabeAnalysisofWiener1987} and %
Yoshida~\cite{Yoshida1992StatisticsSmallDiffusion,Yoshida1992MLE}, the
following result implies
uniform convergence of our asymptotic expansion of transition density
jointly in the whole state space $E$ for the forward variable $x$, the whole
set $\Theta$ for the parameter $\theta$, and an arbitrary compact
subset $%
K\subset$ $E$ for the backward variable $x_{0}$.

\begin{theorem}
\label{propositionconvergencedensity} Under the Assumptions \ref
{assumptionpositivedefinite} and~\ref{aussumptionboundedderivatives},
the transition density expansion (\ref{PXJ}) satisfies
%
\begin{equation}\label{densityapproximationerrorboundequation}
\sup_{(x,x_{0},\theta)\in E\times K\times\Theta}\bigl\llvert
p_{X}^{(J)}(
\Delta,x|x_{0};\theta)-p_{X}(\Delta,x|x_{0};
\theta)\bigr\rrvert=\mathcal{O}\bigl(\Delta^{({J+1-m})/{2}}\bigr)\hspace*{-25pt}
\end{equation}
as $\Delta\rightarrow0$ for $J\geq m$.
\end{theorem}

\begin{pf}
See Appendix~\ref{appendixproofs}.
\end{pf}

It deserves to note that (\ref{densityapproximationerrorboundequation})
gives a theoretical (not necessarily tight) upper bound estimate of the
uniform approximation error of $p_{X}^{(J)}(\Delta,x|x_{0};\theta
)-p_{X}(\Delta,x|x_{0};\theta)$. The effects of dimensionality can be seen
as resulting from the multiplier $\Delta^{-{m}/{2}}$ in the
expansion %
(\ref{PXJ}), which leads to the error magnitude $\Delta^{({J+1-m})/{2}}$.
When $J$ is taken sufficiently large as $J\geq m$, the uniform error is
controlled by taking $\Delta\rightarrow0$.

\section{Numerical performance of density approximation}
\label{sectionimplementationexamples}

In this section, we employ three representative and analytically tractable
examples (the mean-reverting Ornstein--Uhlenbeck process, the Feller square
root process and the double mean-reverting Ornstein--Uhlenbeck process) with
explicitly known transition densities to demonstrate the numerical
performance of the transition density asymptotic expansion proposed in
Section~\ref{sectionexpansionframework}. For all of the examples
investigated in this and the subsequent sections, we provide the first
several expansion terms calculated from the general formula (\ref
{leadingordergeneral}) and (\ref{formulageneralOmegak}) in the
supplementary material~\cite{autokey73}. Higher order correction terms involved in the
numerical implementation are documented in the form of Mathematica notebook,
which will be provided upon request. The density expansions will be
used in
Monte Carlo analysis for approximate maximum likelihood estimation in
Section~\ref{SectionAMLEMC}.

The mean-reverting Ornstein--Uhlenbeck process (also known as the Vasicek
model in financial applications) labeled as MROU is specified as:

\begin{model}
\label{modelvasicek} The MROU (mean-reverting Ornstein--Uhlenbeck) model,
\[
dX(t)=\kappa\bigl(\alpha-X(t)\bigr)\,dt+\sigma \,dW(t).
\]
\end{model}

The Gaussian nature of the MROU model renders a closed-form
transition density, which serves as a benchmark for explicit comparison with
our asymptotic expansion approximations. In the numerical experiments, we
choose a parameter set $\kappa=0.5$, $\alpha=0.06$ and $\sigma=0.03$
similar to those employed in A{\"{\i}}t-Sahalia \cite
{AitSahalia2002Econometrica}.

The Feller square root process (also known as the Cox--Ingersoll--Ross model
in financial applications) labeled as SQR is specified as:

\begin{model}
\label{modelCIR} The SQR (Feller's square root) model,
\[
dX(t)=\kappa\bigl(\alpha-X(t)\bigr)\,dt+\sigma\sqrt{X(t)}\,dW(t).
\]
\end{model}

The combination of the mean-reverting feature and the
Bessel nature (see, e.g., Chapter XI in Revuz and Yor \cite
{RevuzYor}) renders
closed-form transition densities. In particular, we concentrate on
the case where zero is unattainable, that is, the Feller condition
$2\kappa\alpha-\sigma^{2}>0$ holds; see
Feller~\cite{FellerAnnMath51}. In the numerical experiments, we choose a
parameter set $\kappa=0.5$, $\alpha=0.06$ and $\sigma=0.15$
similar to those employed in A{\"{\i}}t-Sahalia \cite
{AitSahalia2002Econometrica}.

We recall that, for the one-dimensional diffusions investigated in %
A{\"{i}}t-Sahalia~\cite{AitSahalia1999JF,AitSahalia2002Econometrica}
and the so-called reducible multivariate diffusions discussed in
A{\"{\i}}t-Sahalia~\cite{AitSahalia2008AS}, the density expansions
proposed in A{\"{i}}t-Sahalia
\cite{AitSahalia1999JF,AitSahalia2002Econometrica,AitSahalia2008AS}
begin with a so-called \textit{Lamperti} transform, which transforms
the
marginal distribution to locally normal. Whenever applied, let $\gamma
(\cdot;\theta)$ denote such a transform, and let $%
Z(t)=\gamma(X(t);\theta)$ denote the process after the transform.
Thus, taking one-dimensional cases as an example, the expansion for
the transition density of $X$ can be constructed from
%
\begin{equation}\label{CIRcombinationmethod}
p_{X}^{(J)}(\Delta,x|x_{0};\theta):=\sigma(x;
\theta)^{-1}p_{Z}^{(J)}\bigl(\Delta,\gamma(x;
\theta)|\gamma(x_{0};\theta);\theta\bigr).
\end{equation}
As momentarily demonstrated in the numerical results, a combination of
the \textit{Lamperti} transform and our expansion leads to faster
convergence, compared with the direct expansion. A heuristic reason for
this phenomenon is as follows. As seen from Section \ref
{sectionexpansionframework}, our expansion is carried out around a
normal distribution. After a \textit{Lamperti} transform, the diffusion
behaves locally as a Brownian motion, which facilitates the
convergence. Therefore, the \textit{Lamperti} transform may accelerate
the convergence of expansion, and thus it is recommended to apply it
whenever it exists and is explicit.
%
%

For multivariate cases, we employ a popular double mean-reverting
Orn\-stein--Uhlenbeck model (see, e.g., A{\"{\i}}t-Sahalia \cite
{AitSahalia2008AS}) labeled as
DMROU, whose transition density is bivariate correlated normal:

\begin{model}
\label{modelA02} The DMROU (double mean-reverting
Ornstein--Uhlenbeck)
mod\-el,
\[
d \pmatrix{ X_{1}(t)
\cr
X_{2}(t)} = \pmatrix{
\kappa_{11} & 0
\cr
\kappa_{21} & \kappa_{22} }
\biggl( \pmatrix{ \alpha_{1}
\cr
\alpha_{2}} -\pmatrix{
X_{1}(t)
\cr
X_{2}(t)} \biggr) \,dt+d \pmatrix{
W_{1}(t)
\cr
W_{2}(t)},
\]
where $\{(W_{1}(t),W_{2}(t))\}$ is a standard two-dimensional Brownian
motion.
\end{model}

According to the classification in
Dai and Singleton~\cite{DaiSingleton1997}, the DMROU model is a
multivariate affine
diffusion process of the $A_{0}(2)$ type. In the numerical
experiments, we choose the parameters as $\kappa_{11}=5$, $\kappa
_{21}=1$, $%
\kappa_{22}=10$ and $\alpha_{1}=\alpha_{2}=0$ similar to those
employed
in A{\"{\i}}t-Sahalia~\cite{AitSahalia2008AS}.

\begin{figure}

\includegraphics{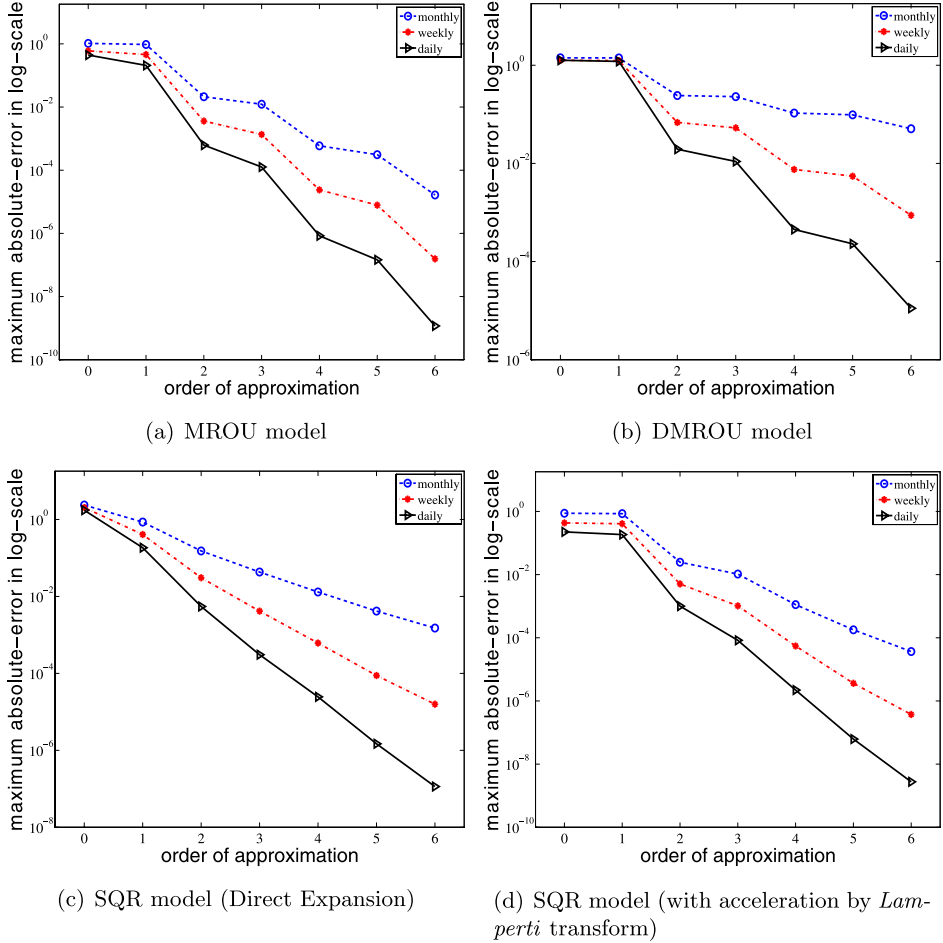}

\caption{Maximum absolute errors of density approximation for Models
\protect\ref{modelvasicek}, \protect\ref{modelCIR} and
\protect\ref{modelA02}.}
\label{fig1}
\end{figure}

Based on the explicit expressions of the true transition densities of the
three models (see, e.g.,
A{\"{i}}t-Sahalia\vadjust{\goodbreak} \cite
{AitSahalia1999JF,AitSahalia2002Econometrica,AitSahalia2008AS}), we
exhibit the error of
$J$th order approximation $e_{X}^{(J)}(\Delta,x|x_{0};\theta
)=p_{X}^{(J)}(\Delta,x|x_{0};\theta)-p_{X}(\Delta,x|x_{0};\theta)$ for
the time increment $\Delta$. The numerical investigation is performed
at a
region $\mathcal{D}$, which is several standard deviations around the mean
of the forward position (i.e., $\mathbb{E} ( X(\Delta
)|X(0)=x_{0} )
$), and as an indicator of the overall performance, the uniform error $%
\max_{x\in\mathcal{D}}|e_{X}^{(J)}(\Delta,x|x_{0};\theta)|$ is
considered. In Figure~\ref{fig1}(a),
(b), (c) and~(d), the uniform errors for the above three
%
benchmark models (MROU, SQR and \mbox{DMROU}) are plotted for monthly, weekly
and daily
monitoring frequencies ($\Delta=1/12,1/52,1/252$) and different orders of
approximation ($J=1,2,3,\ldots,6$). Especially for the SQR model,
the plots are provided for both a direct expansion in Figure
\ref{fig1}(c) and an expansion with \textit{Lamperti}
transform acceleration [see (\ref{CIRcombinationmethod})] in Figure
\ref{fig1}(d).
Such numerical evidence demonstrates that the approximation error
tends to
decrease as the monitoring increment shrinks ($\Delta$ decreases) or more
correction terms are included ($J$ increases), and that the combination
with \textit{Lamperti} transform may accelerate the convergence. As seen
from the dynamics of the SQR model, the volatility function $\sigma
(x)=\sigma\sqrt{x}$ violates Assumption \ref
{aussumptionboundedderivatives} at the point $x=0$. However, the numerical
performance exhibited in Figure~\ref{fig1}(c) and
(d) suggests that the technical assumptions
given in
Section~\ref{SectionModelMLE} are sufficient but not necessary in
order to
guarantee numerical convergence of the density expansion and the resulting
properties of the approximate MLE. From theoretical perspectives, the
singularity at $x=0$ may lead to a significant challenge in mathematically
verifying the convergence of transition density expansion, which can be
regarded as a future research topic.

In the supplementary material~\cite{autokey73}, we document detailed
performance of the
density approximation for the MROU, SQR (for both the direct expansion and
the accelerated approach via \textit{Lamperti} transform) and DMROU models,
respectively. For the former two one-dimensional cases, that is, the
MROU and
SQR models, we plot the errors of approximation corresponding to weekly
monitoring frequency and orders ranging from $J=1,2,\ldots,6$. For the
latter set of graphs, we plot the contours of the approximation errors for
the DMROU model corresponding to weekly monitoring frequency and orders
ranging from $J=1,2,\ldots,6$.

The asymptotic expansion proposed in this paper is essentially
different from that in A{\"{i}}t-Sahalia
\cite{AitSahalia1999JF,AitSahalia2002Econometrica,AitSahalia2008AS} and
other existing large-deviations-based results. First, the expansion
proposed here includes correction terms corresponding to any order of
$\epsilon=\sqrt{\Delta}$; however, in A{\"{i}}t-Sahalia
\cite{AitSahalia1999JF,AitSahalia2002Econometrica,AitSahalia2008AS} and
other methods, expansions include only integer orders of $\Delta$
(even orders of $\epsilon=\sqrt{\Delta}$). Second, the expansion
terms in A{\"{i}}t-Sahalia
\cite{AitSahalia1999JF,AitSahalia2002Econometrica,AitSahalia2008AS}
appear to be longer than the corresponding orders in the expansion
proposed in this paper. Taking the MROU and the SQR model, for example,
the mean-reverting correction starts from the leading order in
A\"{\i}t-Sahalia's expansion; however, in our expansion, the leading
order term is the density of a normal distribution, and the first
appearance of mean-reverting drift parameters is deferred to the
correction term corresponding to first order of $\epsilon$.

%
\begin{figure}

\includegraphics{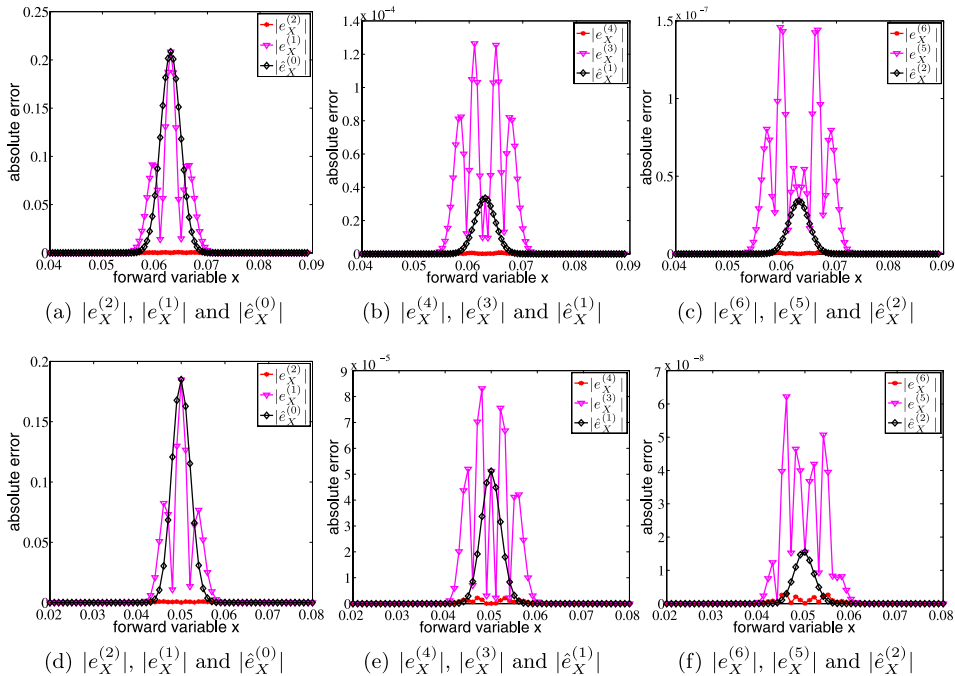}

\caption{Cross comparisons of absolute approximation errors (corresponding
to different orders of expansion) of the proposed method ($|e_{X}^{(J)}|$)
and that of A{\"{\i}}t-Sahalia~\cite{AitSahalia2002Econometrica}
($|\widehat{e}{}^{(J)}_{X}|$).}\label{fig2}
\end{figure}

Let $\widehat{e}{}^{(J)}_{X}(\Delta,x|x_{0};\theta)=\widehat{p}{}^{(J)}
_{X}(\Delta,x|x_{0};\theta)-p_{X}(\Delta,x|x_{0};\theta)$ denote
the approximation error of A\"{\i}t-Sahalia's $J$th order expansion,
where $%
\widehat{p}{}^{(J)}_{X}$ is defined in equation (2.14) as in %
A{\"{\i}}t-Sahalia~\cite{AitSahalia2002Econometrica}. For the method
proposed in this paper, approximation errors are denoted by
$e_{X}^{(J)}(\Delta,x|x_{0};\theta)=p_{X}^{(J)}(\Delta,x|x_{0};\theta)-p_{X}(\Delta,x|x_{0};\theta)$.
Without loss of generality, I employ the MROU and the SQR models to
numerically illustrate the comparison of errors resulting from the
method of A{\"{i}}t-Sahalia
\cite{AitSahalia1999JF,AitSahalia2002Econometrica,AitSahalia2008AS} and
those from this paper. Considering different expressions and
arrangements of correction terms, I make a cross comparison of absolute
errors for different orders from the two methods as exhibited in Figure
\ref{fig2}(a), (b) and (c) for the MROU model as well as Figure
\ref{fig2}(d), (e) and (f) for the SQR model. In particular, we
consider the \textit{Lamperti} transform acceleration for the SQR model
in order to parallel the method in A{\"{\i}}t-Sahalia \cite
{AitSahalia2002Econometrica}. In the comparison, the orders range from
$J=1,2,3,\ldots,6$ for our method,
while $%
J=0,1,2$ for that of A{\"{\i}}t-Sahalia~\cite
{AitSahalia2002Econometrica}. Without loss of
generality, the monitoring frequency is chosen as $\Delta=1/52$. As we will
see, the absolute errors resulting from each two consecutive orders $J=2K-1$
and $J=2K$ of the expansion proposed in this paper sandwich that resulting
from the order $K-1$ of the expansion proposed in %
A{\"{\i}}t-Sahalia~\cite{AitSahalia2002Econometrica}, for $K=0,1,2$.
The two methods both
admit small magnitude of errors resulting from low order approximations and
are comparable to each other as more correction terms are included.

\section{Approximate maximum-likelihood estimation}
\label{SectionAMLE}

This section is devoted to a method of approximate MLE based on the
asymptotic expansion for transition density proposed in Section \ref
{sectionexpansionframework}. Similar to A{\"{i}}t-Sahalia
\cite{AitSahalia1999JF,AitSahalia2002Econometrica,AitSahalia2008AS},
the $J$th order expansion of the log-density can be given by
\[
l_{X}^{(J)}(\Delta,x|x_{0};\theta):=-
\frac{m}{2}\log\Delta+\log\bigl[ \det D(x_{0}) \bigr] +\sum
_{k=0}^{J}\Lambda_{k} \biggl(
D(x_{0})\frac
{x-x_{0}%
}{\sqrt{\Delta}} \biggr) \epsilon^{k}
\]
for any $J=0,1,2,\ldots\,$, where the correction terms $\Lambda_{k}$
can be explicitly calculated from straightforward differentiation of
the density expansion (\ref{PXJ}).

\begin{figure}

\includegraphics{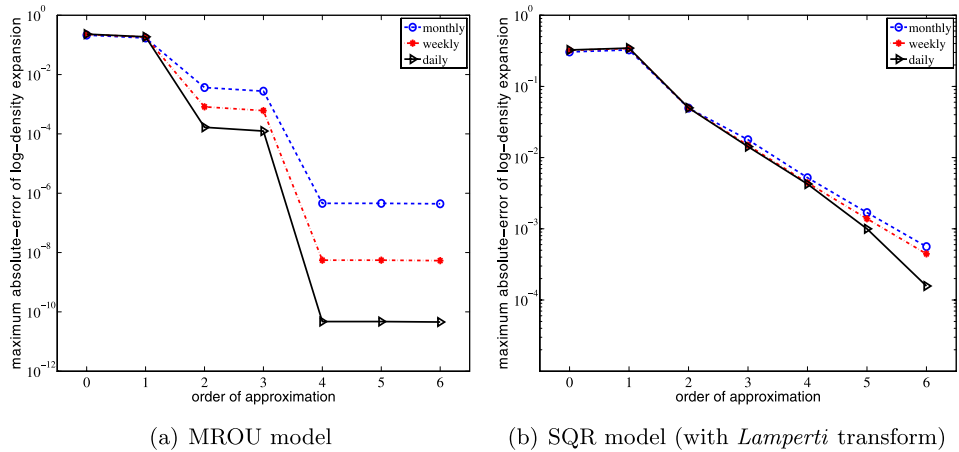}

\caption{Maximum absolute errors of log-density approximation for
Models \protect\ref{modelvasicek} and \protect\ref{modelCIR}.}
\label{fig3}
\end{figure}

Without loss of generality, we employ Model~\ref{modelvasicek} (MROU)
and Model~\ref{modelCIR} (SQR) to illustrate the convergence of uniform
errors of the log-density expansions ($\max_{x\in\mathcal
{D}}|l_{X}^{(J)}(\Delta,x|x_{0};\theta)-\log
p_{X}(\Delta,x|x_{0};\theta)|$) in Figure~\ref{fig3}(a)
and (b) in a similar way as Figure
\ref{fig1}(a)--(d) do for the
uniform errors of density expansions. For the MROU model, Figure
\ref{fig3}(a) shows the uniform errors of its
log-density expansions. For the SQR model, Figure
\ref{fig3}(b) plots the uniform errors of its
\textit{Lamperti}-transformed log-density expansions which are
naturally calculated from
%
\begin{equation}\label{loglikelihoodexpansion}
l_{X}^{(J)}(\Delta,x|x_{0};\theta):=-\log
\sigma(x;\theta)+l_{Z}^{(J)}\bigl(\Delta,\gamma(x;\theta)|
\gamma(x_{0};\theta);\theta\bigr),
\end{equation}
where $\gamma$ is the \textit{Lamperti} transform and $Z(t)=\gamma
(X(t);\theta)$.

%

By analogy to the log-likelihood function $\ell_{n}(\theta)$ in %
(\ref{loglikelihood}), we introduce the $J$th order approximate
log-likelihood function
%
\begin{equation}\label{approximateloglikelihoodfunction}
\ell_{n}^{(J)}(\theta)=\sum_{i=1}^{n}L_{i}^{(J)}(
\theta),
\end{equation}
where $L_{i}^{(J)}(\theta)=l_{X}^{(J)}(\Delta,X(i\Delta
)|X((i-1)\Delta
);\theta)$. According to Assumption \ref
{assumption3timesdifferentiable}, we assume, for simplicity, that the
true log-likelihood function $%
\ell_{n}(\theta)$ admits a unique maximizer $\widehat{\theta
}_{n}$, which
serves as the true maximum-likelihood estimator. Similarly, let
$\widehat{%
\theta}{}^{(J)}_{n}$ be the approximate maximum-likelihood estimator
of order $J$ obtained from maximizing $\ell_{n}^{(J)}(\theta)$.
Setting up or refining technical conditions for ensuring the
identification and the usual asymptotic properties of the true but
generally incomputable MLE is beyond the scope of this paper, and
can be investigated as a future research topic; see, for example,
A{\"{\i}}t-Sahalia~\cite{AitSahalia2002Econometrica} for the
discussion of
one-dimensional cases. As a consequence of Theorem \ref
{propositionconvergencedensity}, we set up the convergence of the
approximate MLE $\widehat{\theta}{}^{(J)}_{n}$ to the true MLE
$\widehat{%
\theta}_{n}$ in what follows.

\begin{proposition}
\label{theoremconvergenceAMLE}Under Assumptions \ref
{assumptionpositivedefinite} and~\ref{aussumptionboundedderivatives},
the approximate maximum-likelihood estimator obtained from optimizing
(\ref{approximateloglikelihoodfunction}) satisfies that, for the fixed sample
size $n$,
%
\begin{equation}\label{AMLEtoMLE}
\widehat{\theta}{}^{(J)}_{n}-\widehat{\theta}_{n}
\stackrel{P} {\rightarrow}0
\end{equation}
as $\Delta\rightarrow0$ for $J\geq m$.
\end{proposition}

\begin{pf}
See Appendix~\ref{appendixproofs}.
\end{pf}

Though the convergence in (\ref{AMLEtoMLE}) is theoretically justified as
the monitoring increment $\Delta$ shrinks to $0$ for any fixed order $J$,
the convergence of (\ref{AMLEtoMLE}) may also hold as $J\rightarrow
\infty$
for a range of fixed values of $\Delta$. This is analogous to the Taylor
expansion in classical calculus. The respective effects on the discrepancy
between the approximate MLE and the true MLE resulting from shrinking $%
\Delta$ and increasing expansion orders $J$ are illustrated via Monte Carlo
evidence in Section~\ref{SectionAMLEMC}. In particular, we will
demonstrate numerically that for an arbitrary $\Delta$, a larger order $J$
results in a better approximation of the MLE.


\section{Monte Carlo evidence of approximate maximum-likelihood estimation}
\label{SectionAMLEMC}

To further demonstrate the convergence issues
discussed in the previous sections, we provide Monte Carlo evidence of
approximate maximum-likelihood estimation for the three models
discussed in Section~\ref{sectionimplementationexamples}. Let $N$
denote the number of sample paths generated from the transition
distributions; let $n$ denote the number of observations on each path.
For finite-sample results, we report the mean and standard deviation of
the discrepancy between the MLE and the true parameter value (i.e.,
$\widehat{\theta}_{n}-\theta^{\mathrm{True}}$), and the discrepancy
between the approximate MLE and the MLE (i.e., $\widehat{\theta
}{}^{(J)}_{n}-\widehat{\theta}_{n}$).

For the two Gaussian models, that is, the MROU model and the DMROU
model, the
situation considered here is restricted to the stationary case. Therefore,
the asymptotic variance of the maximum-likelihood estimator is given by the
inverse of Fisher's information matrix, which is the lowest possible
variance of all estimators. So, as in
A{\"{i}}t-Sahalia \cite
{AitSahalia1999JF,AitSahalia2002Econometrica,AitSahalia2008AS}, we
assume that $\kappa>0$ for the MROU model, and $%
\kappa_{11}>0$ and $\kappa_{22}>0$ for the DMROU model. By the nature of
stationarity, one has the local asymptotic normal structure for the
maximum-likelihood estimator~$\widehat{\theta}_{n}$, that is,
%
\begin{equation}\label{asymptoticdistribution}
\sqrt{n}(\widehat{\theta}_{n}-\theta)\stackrel{\mathcal{D}} {
\rightarrow}%
\mathcal{N}\bigl(0,i(\theta)^{-1}\bigr)
\end{equation}
as $n\rightarrow\infty$ with $\Delta$ fixed. Here, the Fisher information
matrix is calculated~as
%
\begin{equation}\label{Fisherinformationmatrix}
i(\theta)=-\mathbb{E} \biggl( \frac{\partial^{2}L_{1}(\theta
)}{\partial
\theta\,\partial\theta^{T}} \biggr),
\end{equation}
where $^{T}$ denotes matrix transposition.

%
\begin{table}
\tabcolsep=0pt
\caption{Monte Carlo evidence for the MROU model}
\label{tableVasicek}
\begin{tabular*}{\tablewidth}{@{\extracolsep{\fill}}lccccccd{2.6}c@{}}
\hline
& \multicolumn{2}{c}{\textbf{Asymptotic}} &
\multicolumn{2}{c}{\textbf{Finite sample}} &
\multicolumn{2}{c}{\textbf{Finite sample}} &
\multicolumn{2}{c@{}}{\textbf{Finite sample}} \\
\multirow{2}{45.5pt}[-5.7pt]{\textbf{Parameters}
$\bolds{\theta^{\mathrm{True}}}$} & \multicolumn{2}{c}{$\bolds{\widehat{\theta}_{n}
-\theta^{\mathrm{True}}}$} & \multicolumn{2}{c}{$\bolds{\widehat{\theta}_{n}
-\theta^{\mathrm{True}}}$} & \multicolumn{2}{c}{$\bolds{\widehat{\theta
}{}^{(3)}_{n}-\widehat{\theta}%
_{n}}$} & \multicolumn{2}{c@{}}{$\bolds{\widehat{\theta}{}^{(6)}_{n}-\widehat
{\theta}%
_{n}}$} \\[-4pt]
& \multicolumn{2}{c}{\hrulefill} & \multicolumn{2}{c}{\hrulefill}
& \multicolumn{2}{c}{\hrulefill} & \multicolumn{2}{c@{}}{\hrulefill}\\
& \textbf{Mean} & \textbf{Stddev} &
\textbf{Mean} & \textbf{Stddev} & \textbf{Mean} & \multicolumn{1}{c}{\textbf{Stddev}} & \multicolumn{1}{c}{\textbf{Mean}}
& \multicolumn{1}{c@{}}{\textbf{Stddev}} \\
\hline
$\Delta=1/52$ & & & & & & & & \\
\quad$\kappa=0.5$ & 0 & 0.229136 & 0.245175 & 0.329396 &
0.013477 & 0.014645 & 0.000002 & 0.000102 \\
\quad$\alpha=0.06$ & 0 & 0.013682 & 0.000329 &
0.015202 &
0.000002 & 0.000318 & 0.000000 & 0.000003 \\
\quad$\sigma=0.03$ & 0 & 0.000674 & 0.000021 &
0.000675 &
0.000003 & 0.000015 & -0.000000 & 0.000000 \\
[4pt]
$\Delta=1/12$ & & & & & & & & \\
\quad$\kappa=0.5$ & 0 & 0.111867 & 0.054162 & 0.124773 &
0.028923 & 0.014382 & -0.000003 & 0.000297 \\
\quad$\alpha=0.06$ & 0 & 0.006573 & 0.000097 &
0.006440 &
0.000002 & 0.000174 & 0.000000 & 0.000014 \\
\quad$\sigma=0.03$ & 0 & 0.000685 & 0.000022 &
0.000687 &
0.000025 & 0.000022 & 0.000000 & 0.000001 \\
\hline
\end{tabular*}
\legend{Notes. The number of
simulation trials is $N=5000$ and the number of observations on each
path is $n=1000$.}
\end{table}

Without loss of generality, we analyze the results of the MROU
model in what follows. As seen from Table~\ref{tableVasicek}, the
asymptotic distribution of $\widehat{\theta}%
_{n}-\theta^{\mathrm{True}}$ is calculated from %
(\ref{asymptoticdistribution}) and
(\ref{Fisherinformationmatrix}). The small discrepancy between the
finite-sample and asymptotic standard deviations of $\widehat{\theta
}_{n}-\theta^{\mathrm{True}}$ indicates that the choice of sample
size $n=1000$ is approaching an optimality. When the monitoring
frequency $\Delta$ shrinks, or when the order of approximation $J$
increases, the approximate
MLEs obtained from maximizing the approximate log-likelihood function %
(\ref{approximateloglikelihoodfunction}) get closer to the exact (but
usually incomputable) MLEs, and thus get closer to the true parameter,
if the sample size $n$ is large enough. This can be seen by a
comparison of some outputs with relatively larger bias and standard
deviations resulting from relatively lower order expansions or larger
monitoring increments with those improved outputs resulting from
relatively higher order expansions and smaller monitoring increments.
This phenomenon reconciles our discussions in Section~\ref{SectionAMLE}.

While holding the length of sampling interval $\Delta$ fixed, the
approximation error $\widehat{\theta}{}^{(J)}_{n}-\widehat{\theta
}_{n}$
decreases and is dominated by the intrinsic sampling error $\widehat
{\theta}%
_{n}-\theta^{\mathrm{True}}$ as $J$ increases. Therefore, according to
A{\"{i}}t-Sahalia \cite
{AitSahalia1999JF,AitSahalia2002Econometrica,AitSahalia2008AS}, a
small-order approximation (e.g., the $%
\widehat{\theta}{}^{(6)}_{n}$ for the MROU model) is adequate enough for
replacing the true MLE $\widehat{\theta}_{n}$ for the purpose of
estimating unknown parameter $\theta$. According to A{\"{i}}t-Sahalia
\cite{AitSahalia1999JF,AitSahalia2002Econometrica,AitSahalia2008AS},
once the approximation error resulting from replacing the true density
$p_{X}$ by its approximation, say $p_{X}^{(J)}$, is dominated by the
sampling error [usually estimated from asymptotic variance computed
from (\ref{asymptoticdistribution})] due to the true maximum-likelihood
estimation, such $p_{X}^{(J)}$ is appropriate in practice. Such a
proper replacement has an effect that is statistically
indiscernible from the\vspace*{2pt} sampling variation of the true yet incomputable
MLE $%
\widehat{\theta}_{n}$ around $\theta$. As a result of the fast
development of modern computation and optimization technology,
calculation of high-order likelihood approximations will become
increasingly feasible; thus errors between approximate MLE and MLE
can be improved to become arbitrarily small, at least in principle.

Owing to the limited space in this paper and the similarity in the
pattern of results to those of the MROU model, we collect the
simulation results for the DMROU and the SQR models in the
supplementary material~\cite{autokey73}. In particular, for the SQR
model, the
simulation results will demonstrate that a combination of the
\textit{Lamperti} transform and our expansion may enhance the
efficiency of the estimation. Moreover, in the
supplementary material~\cite{autokey73}, we will investigate two more
sophisticated data-generating processes (arising from financial
modeling) with rich drift and diffusion specifications, in which the
Lamperti transform either requires computationally demanding
implicit integration and inversion or does not exist due to a
multivariate irreducible specification; see
A{\"{\i}}t-Sahalia~\cite{AitSahalia2008AS}.

\section{Concluding remarks}
\label{sectionconcludingremarks}

This paper contributes a method for approximate maximum-likelihood
estimation (MLE) of multivariate diffusion processes from discretely
sampled data, based on a closed-form asymptotic expansion for
transition density, for which any arbitrary order of corrections can
be systematically obtained through a generally implementable
algorithm. Numerical examples and Monte Carlo evidence for
illustrating the performance of density asymptotic expansion and the
resulting approximate MLE are provided in order to demonstrate the
wide applicability of the method. Based on some sufficient (but not
necessary) technical conditions, the convergence and asymptotic
properties are theoretically justified. Owing to the limited space
of this paper which focuses on introducing a method of estimation,
investigations on more asymptotic properties related to the
approximate MLE can be regarded as a future
research topic, for example, a tighter upper bound for the discrepancy (\ref
{AMLEtoMLE}) based on the error estimate of the transition density expansion
(\ref{densityapproximationerrorboundequation}), as well as the
consistency and asymptotic distribution of the approximate MLE under various
sampling schemes in terms of monitoring frequency and observational horizon;
see, for example, Yoshida~\cite{Yoshida1992Discrete}, Kessler \cite
{Kessler1997}, %
Bibby and S{\o}rensen~\cite{BibbySorensen1995} and Genon-Catalot and
Jacod~\cite{GenonCatalotJacod1993}. In this
regard, we note Chang and Chen~\cite{ChangChen2011} for analyzing the
asymptotic
properties of the approximate MLE proposed in %
A{\"{\i}}t-Sahalia~\cite{AitSahalia2002Econometrica} of
one-dimensional diffusion
processes. One may also apply the idea for explicitly approximating
transition density in various other aspects of statistical inference, for
which explicit asymptotic expansions of certain quantities are helpful.


\begin{appendix}\label{app}
\section{\texorpdfstring{Explicit calculation of conditional expectation (\lowercase{\protect\ref{ce_product_stratonovichintegrals}})}
{Explicit calculation of conditional expectation (3.24)}}
\label{subsectionCE}

In this section, we expatiate on a general algorithm for explicitly
calculating the conditional expectation
(\ref{ce_product_stratonovichintegrals}) of multiplication of
iterated
Stratono\-vich integrals as a multivariate polynomial in $%
z=(z_{1},z_{2},\ldots,z_{m})\in\mathbb{R}^{m}$. In addition to
theoretical interests, iterated stochastic integral plays important
roles in many applications arising from stochastic modeling, for example,
the analysis of convergence rate of various methods for
approximating solutions to stochastic differential equations; see
Kloeden and Platen~\cite{KloedenPlaten1999}. Special cases for conditional
expectations of iterated It\^{o} stochastic integrals (without
integral with respect to the time variable) can be found in, for example,
Nualart, {\"U}st{\"u}nel and Zakai~\cite{NualartUstunelZakai1988},
Yoshida~\cite{Yoshida1992StatisticsSmallDiffusion,Yoshida1992MLE} and
Kunitomo and Takahashi~\cite{KunitomoTakahashi2001mathfinance}.

To present our algorithm, similar to the definition of iterated Stratonovich
integral, we define
%
\begin{equation}\label{ItoIntfgeneral}\quad
I_{\mathbf{i}}[f](t):=\int_{0}^{t}\int
_{0}^{t_{1}}\cdots\int_{0}^{t_{n-1}}f(t_{n})\,dW_{i_{n}}(t_{n})\cdots dW_{i_{2}}(t_{2})\,dW_{i_{1}}(t_{1})
\end{equation}
as an iterated It\^{o} integral for an arbitrary index $\mathbf{i}%
=(i_{1},i_{2},\break\ldots,i_{n})\in\{0,1,2,\ldots, m\}^{n}$ with a
right-continuous integrand $f$. To lighten the notation, the integral
$I_{%
\mathbf{i}}[1](t)$ is abbreviated to $I_{\mathbf{i}}(t)$.

Before discussing details in the following subsections, we briefly
outline a general algorithm, which can be implemented using any
symbolic packages, for example, Mathematica. Throughout our discussion, the
iterated (Stratonovich or It\^{o}) stochastic integrals may involve
integrations with respect to not only Brownian motions but also time
variables.

\begin{algorithm*}
\label{algorithmCE}
\begin{itemize}
\item Convert each iterated Stratonovich integral in
(\ref{ce_product_stratonovichintegrals}) to a linear combination of
iterated It\^{o} integrals;

\item Convert each multiplication of iterated It\^{o} integrals resulting
from the previous step to a linear combination of iterated It\^{o} integrals;

\item Compute the conditional expectation of iterated It\^{o} integral via
an explicit construction of Brownian bridge.
\end{itemize}
\end{algorithm*}
%
\subsection{Conversion from iterated Stratonovich integrals to It\^{o} integrals}

Denote by $l(\mathbf{i}%
):=l((i_{1},\ldots,i_{n}))=n$ the length of the index $\mathbf{i}$.
Denote by $-\mathbf{i}$ an index obtained from deleting the first
element of $\mathbf{i}$. In particular, if\vadjust{\goodbreak}
$l(\mathbf{i})=0$, we define $J_{%
\mathbf{i}}[f](t)=f(t)$ by slightly extending the definition
(\ref{StratIntfgeneral}). According to page 172 of
Kloeden and Platen~\cite{KloedenPlaten1999}, we have the following conversion
algorithm: for the case of $l(\mathbf{i}%
)=0$ or $1$, we have $J_{\mathbf{i}}(t)=I_{\mathbf{i}}(t)$; for the
case of $l(\mathbf{i})\geq2$, we have
%
\begin{equation}\label{conversionSTIT}
J_{\mathbf{i}}(t)=I_{(i_{1})}\bigl[J_{-\mathbf{i}}(\cdot)
\bigr](t)+1_{ \{
i_{1}=i_{2}\neq0 \} }I_{(0)} \bigl[ \tfrac{1}{2}J_{-(-\mathbf
{i})}(
\cdot) \bigr] (t).
\end{equation}
For example, if $l(\mathbf{i})=2$, one has
\[
{J_{\mathbf{i}}(t)=I_{_{\mathbf{i}}}(t)+\tfrac{1}{2}1}_{%
{\{i_{1}=i_{2}\neq0\}}}{
I_{(0)}(t).}
\]
Thus, with the conversion algorithm (\ref{conversionSTIT}), we convert
each iterated Strato\-novich integral in (\ref
{ce_product_stratonovichintegrals}) to a linear combination of iterated
It\^{o} integrals. Thus, the product $\prod_{\omega
=1}^{l}J_{\mathbf{%
i}_{\omega}}(1)$ can be expanded as a linear combination of multiplication
of It\^{o} integrals.

\subsection{Conversion from multiplication of It\^{o} integrals to a linear
combination}
\label{subsectionconversionmultiplicationtolinearcombination}

We provide a simple recursion algorithm for converting a multiplication
of iterated It\^{o} integrals to a linear combination. According to
Lemma~2 in Tocino~\cite{Tocino2009}, a product of two It\^{o} integrals
as defined in (\ref{ItoIntfgeneral}) satisfies that
%
\begin{eqnarray}\label{conversionProdJtolinearcombinationJ}
I_{\bolds{\alpha}}(t)I_{\bolds{\beta}}(t)&=&\int_{0}^{t}I_{\bolds
{\alpha}%
}(s)I_{-\bolds{\beta}}(s)\,dW_{\beta_{1}}(s)+
\int_{0}^{t}I_{-\bolds
{\alpha
}}(s)I_{\bolds{\beta}}(s)\,dW_{\alpha_{1}}(s)
\nonumber\\[-8pt]\\[-8pt]
&&{}+\int_{0}^{t}I_{-\bolds{%
\alpha}}(s)I_{-\bolds{\beta}}(s)1_{\{\alpha_{1}=\beta_{1}\neq
0\}}\,ds\nonumber
\end{eqnarray}
for any arbitrary indices $\bolds{\alpha}=(\alpha
_{1},\alpha_{2},\ldots,\alpha_{p})$ and $\bolds{\beta}=(\beta
_{1},\beta_{2},\ldots,\beta_{q})$. Iterative applications of this
relation render a linear combination form of $I_{\bolds{\alpha
}}(t)I_{\bolds{\beta}}(t)$. Inductive applications of such an
algorithm convert a product of any number of iterated It\^{o}
integrals to a linear combination. Therefore, our immediate task is
reduced to the calculation of conditional expectations of iterated
It\^{o} integrals.

\subsection{Conditional expectation of iterated It\^{o} integral}

We focus on the explicit calculation of conditional expectations
of the following type:
%
\begin{eqnarray}\label{CEIto}\qquad
&&\mathbb{E}\bigl(I_{\mathbf{i}}(1)|W(1)=z\bigr)
\nonumber\\[-8pt]\\[-8pt]
&&\qquad=\mathbb{E} \biggl( \int_{0}^{1}\int
_{0}^{t_{1}}\cdots\int_{0}^{t_{n-1}}dW_{i_{n}}({t_{n})}\cdots
dW_{i_{2}}(t_{2})\,dW_{i_{1}}(t_{1})|W(1)=z \biggr). \nonumber
\end{eqnarray}
By an explicit construction of Brownian bridge (see page 358 in %
Karatzas and Shreve~\cite{KaraztasShreve}), we obtain the following
distributional
identity, for any $k=1,2,\ldots,m$:
\[
\bigl( W_{k}(t)|W(1)=z \bigr) \stackrel{\mathcal{D}} {=} \bigl(
W_{k}(t)|W_{k}(1)=z_{k} \bigr) \stackrel{\mathcal
{D}} {=}BB_{k}^{z}(t):=%
\mathcal{B}_{k}(t)-t
\mathcal{B}_{k}(1)+tz_{k},
\]
where $\mathcal{B}_{k}$'s are independent Brownian motions and $%
BB_{k}^{z}(t):=\mathcal{B}_{k}(t)-t\mathcal{B}_{k}(1)+tz_{k}$ is
distributed as a Brownian bridge starting from $0$ and ending at
$z_{k}$ at time~$1$. For ease of exposition, we also introduce
$\mathcal{B}_{0}(t)\equiv0$ and $z_{0}=1$. Therefore, the
condition $W(1)=z$ in (\ref{CEIto}) can be
eliminated since
%
\begin{eqnarray}\label{convertedCE}
&&\mathbb{E} \biggl( \int_{0}^{1}\int
_{0}^{t_{1}}\cdots\int_{0}^{t_{n-1}}dW_{i_{n}}({t_{n})}\cdots
dW_{i_{2}}(t_{2})\,dW_{i_{1}}(t_{1})|W(1)=z \biggr) \nonumber\\
&&\qquad=\mathbb{E} \biggl( \int_{0}^{1}\int
_{0}^{t_{1}}\cdots\int_{0}^{t_{n-1}}d
\bigl(%
\mathcal{B}_{i_{n}}(t_{n})-t_{n}
\mathcal{B}_{i_{n}}(1)+t_{n}z_{i_{n}}\bigr)\cdots\nonumber\\[-8pt]\\[-8pt]
&&\hspace*{121.7pt}d
\bigl(%
\mathcal{B}_{i_{2}}(t_{2})
-t_{2}\mathcal{B}_{i_{2}}(1)+t_{2}z_{i_{2}}
\bigr) \nonumber\\
&&\hspace*{140pt}d\bigl(\mathcal{B}_{i_{1}}(t_{1})-t_{1}
\mathcal{B}_{i_{1}}(1)+t_{1}z_{i_{1}}\bigr)%
\biggr).\nonumber
\end{eqnarray}
An early attempt using the idea of Brownian bridge to deal with
conditional expectation (\ref{CEIto}) can be found in
Uemura~\cite{Uemura1987}, which investigated the calculation of heat
kernel expansion in the diagonal case. It is worth mentioning that,
instead of giving a method for explicitly calculating
(\ref{CEIto}), Uemura~\cite{Uemura1987} employed discretization of
stochastic integrals to show that (\ref{CEIto}) has the structure
of a multivariate polynomial in $z$ with unknown coefficients.
Therefore, the validity of the above derivation can be seen from the
definition of stochastic integral as a limit
of discretized summation. In particular, the random variables $\mathcal
{B}%
_{i_{1}}(1),\mathcal{B}_{i_{2}}(1),\ldots,\mathcal{B}_{i_{n}}(1)$
are not
involved in the integral in (\ref{convertedCE}). The integrals with respect
to $d\mathcal{B}_{i_{k}}(t_{k})$ are in the sense of usual stochastic
integrals; the integrals with respect to $dt_{k}$ are in the sense of
Lebesgue integrals.

By expanding the right-hand side of (\ref{convertedCE}) and collecting
terms according to monomials of $z_{i}$'s, we express (\ref{CEIto}) as a
multivariate polynomial in $z$:
\[
\mathbb{E}\bigl(I_{\mathbf{i}}(1)|W(1)=z\bigr)=\sum
_{k=0}^{n}%
\sum_{\{l_{1},l_{2},\ldots,l_{k}\}\subset\{1,2,\ldots,n\}}
c(l_{1},l_{2},\ldots,l_{k})z_{i_{l_{1}}}z_{i_{l_{2}}}\cdots
z_{i_{l_{k}}},
\]
where the coefficients are determined by
%
\begin{eqnarray}\label{cLk}
&&
c(l_{1},l_{2},\ldots,l_{k})\nonumber\\
&&\qquad:=\mathbb{E}\int
_{0}^{1}\int_{0}^{t_{1}}\cdots\int_{0}^{t_{n-1}}d \bigl(
\mathcal{B}_{i_{n}}(t_{n})-t_{n}
\mathcal{B}%
_{i_{n}}(1) \bigr)\cdots
\nonumber\\
&&\qquad\hspace*{97pt} d\bigl(\mathcal{B}_{i_{l_{k}+1}}(t_{l_{k}+1})-t_{l_{k}+1}
\mathcal{B}%
_{i_{l_{k}+1}}(1)\bigr)\nonumber\\
&&\qquad\hspace*{97pt} dt_{l_{k}}\,d\bigl(
\mathcal{B}%
_{i_{l_{k}-1}}(t_{l_{k}-1})-t_{l_{k}-1}
\mathcal{B}_{i_{l_{k}-1}}(1)\bigr)\cdots
\nonumber\\[-8pt]\\[-8pt]
&&\qquad\hspace*{97pt}  d\bigl(\mathcal{B}_{i_{l_{2}+1}}(t_{l_{2}+1})-t_{l_{2}+1}
\mathcal{B}%
_{i_{l_{2}+1}}(1)\bigr)\nonumber\\
&&\qquad\hspace*{97pt} dt_{l_{2}}\,d\bigl(
\mathcal{B}%
_{i_{l_{2}-1}}(t_{l_{2}-1})-t_{l_{2}-1}
\mathcal{B}_{i_{l_{2}-1}}(1)\bigr)\cdots
\nonumber
\\
&&\qquad\hspace*{97pt}  d\bigl(\mathcal{B}_{i_{l_{1}+1}}(t_{l_{1}+1})-t_{l_{1}+1}
\mathcal{B}%
_{i_{l_{1}+1}}(1)\bigr)\nonumber\\
&&\qquad\hspace*{97pt} dt_{l_{1}}\,d\bigl(
\mathcal{B}%
_{i_{l_{1}-1}}(t_{l_{1}-1})-t_{l_{1}-1}
\mathcal{B}_{i_{l_{1}-1}}(1)\bigr)\cdots
\nonumber
\\
&&\qquad\hspace*{97pt}  d\bigl(\mathcal{B}_{i_{1}}(t_{1})-t_{1}
\mathcal{B}_{i_{1}}(1)\bigr).
\nonumber
\end{eqnarray}
Algebraic calculation from expanding the terms like $d ( \mathcal
{B}%
_{i_{n}}(t_{n})-t_{n}\mathcal{B}_{i_{n}}(1) ) $ simplifies (\ref{cLk})
as a linear combination of expectations of the following form: $\mathbb
{E}(%
\mathcal{B}_{m_{1}}(1)\mathcal{B}_{m_{2}}(1)\cdots\mathcal
{B}_{m_{r}}(1)I_{%
\mathbf{j}}(1))$ where $I_{\mathbf{j}}(1)$ is an iterated It\^{o} integral.

By viewing $\mathcal{B}_{m_{i}}(1)$ as $\int_{0}^{1}d\mathcal{B}%
_{m_{i}}(t_{1})$, we have
%
\begin{equation}\label{expectationproduct}
\mathbb{E}\bigl(\mathcal{B}_{m_{1}}(1)
\mathcal{B}_{m_{2}}(1)\cdots\mathcal{B}_{m_{r}}(1)I_{\mathbf{j}}(1)\bigr)=
\mathbb{E} \Biggl( \prod_{i=1}^{r}\int
_{0}^{1}d\mathcal{B}_{m_{i}}(t_{1})I_{\mathbf{j}%
}(1)
\Biggr).
\end{equation}
To calculate this expectation, we use the algorithm proposed in
Section~\ref{subsectionconversionmultiplicationtolinearcombination} to convert
$%
\prod_{i=1}^{r}\int_{0}^{1}d\mathcal
{B}_{m_{i}}(t_{1})I_{\mathbf{j}%
}(1)$ to a linear combination of iterated It\^{o} integrals.
Finally, we need to calculate expectation of iterated It\^{o}
integrals without conditioning. For any arbitrary index
$\mathbf{i}=(i_{1},i_{2},\ldots,i_{n})\in\{0,1,2,\ldots,m\}^{n}$,
we have
\begin{eqnarray*}
\mathbb{E}I_{\mathbf{i}}(1) &=&\mathbb{E} \biggl( \int
_{0}^{1}\int_{0}^{t_{1}}\cdots\int_{0}^{t_{n-1}}dW_{i_{n}}({t_{n})}\cdots
dW_{i_{2}}(t_{2})\,dW_{i_{1}}(t_{1}) \biggr)
\\
&=&\int_{0}^{1}\int_{0}^{t_{1}}\cdots\int_{0}^{t_{n-1}}d{t_{n}}\cdots
dt_{2}\,dt_{1}\equiv\frac{1}{n!},
\end{eqnarray*}
if $\mathbf{i}=(i_{1},i_{2},\ldots,i_{n})=(0,0,\ldots,0)$ and
$\mathbb{E}%
I_{\mathbf{i}}(1)=0$, otherwise (by the martingale property of
stochastic integrals).
%

\section{Proofs}
\label{appendixproofs}

\subsection{\texorpdfstring{Proof of Theorem \protect\ref{propositiongeneralOmegak}}
{Proof of Theorem 1}}
Using the chain rule and the Taylor theorem, the $k$th ($k\geq1$)
order correction term for $\delta(Y^{\epsilon}-y)$ admits the
following form:
%
\begin{eqnarray}\label{Phikgeneral}
\Phi_{k}(y)&=&\sum_{ ( l,\mathbf{r}(l),\mathbf{j}(l) ) \in
S_{k}}\,
\frac{1}{l!}\partial^{\mathbf{r}}\delta\bigl(B(1)-y\bigr
)Y_{j_{1},r_{1}}Y_{j_{2},r_{2}}\cdots Y_{j_{l},r_{l}}, 
\end{eqnarray}
where $\partial^{\mathbf{r}}$ denotes $\frac{\partial}{\partial
x_{r_{1}}}\,
\frac{\partial}{\partial x_{r_{2}}}\cdots\frac{\partial}{\partial
x_{_{^{r_{l}}}}}$ for simplicity. Thus, taking expectation of (\ref
{Phikgeneral}) and
applying (\ref{Yexpansion}), we obtain that
\begin{eqnarray*}
\Omega_{k}(y) &=&\mathbb{E}\Phi_{k}(y)\\
&=&\sum
_{ ( l,\mathbf
{r}(l),%
\mathbf{j}(l) ) \in S_{k}}\frac{1}{l!}%
D_{r_{1}r_{1}}(x_{0})D_{r_{2}r_{2}}(x_{0})\cdots D_{r_{l}r_{l}}(x_{0})
\\
&&\hspace*{52pt}{}\times\mathbb{E} \bigl( \partial^{\mathbf{r}}\delta\bigl(B(1)-y
\bigr)F_{j_{1}+1,r_{1}}F_{j_{2}+1,r_{2}}\cdots F_{j_{l}+1,r_{l}} \bigr).
\end{eqnarray*}
Employing the integration-by-parts property of the Dirac delta function
(see, e.g., Section 2.6 in Kanwal \cite
{Kanwal2004generalizedfunctions}), the
conditional expectation can be computed as
\begin{eqnarray*}
&&\mathbb{E}\bigl[\partial^{\mathbf{r}}\delta\bigl(B(1)-y
\bigr)F_{j_{1}+1,r_{1}}F_{j_{2}+1,r_{2}}\cdots F_{j_{l}+1,r_{l}}\bigr]
\\
&&\qquad=\int_{b\in\mathbb{R}^{d}}\mathbb{E\bigl[}\partial^{\mathbf
{r}}\delta
\bigl(B(1)-y\bigr)F_{j_{1}+1,r_{1}}F_{j_{2}+1,r_{2}}\cdots
F_{j_{l}+1,r_{l}}|B(1)=b
\bigr]\\
&&\hspace*{23.5pt}\qquad\quad{}\times\phi_{\Sigma(x_{0})}(b)\,db
\\
&&\qquad=(-1)^{l}\,\partial^{\mathbf{r}} \bigl[ \mathbb{E
\bigl[}%
F_{j_{1}+1,r_{1}}F_{j_{2}+1,r_{2}}\cdots F_{j_{l}+1,r_{l}}|W(1)=
\sigma(x_{0})^{-1}D(x_{0})^{-1}y\bigr]\\
&&\qquad\hspace*{256pt}{}\times
\phi_{\Sigma(x_{0})}(y) \bigr],
\end{eqnarray*}
where $\phi_{\Sigma(x_{0})}(y)$ is given in
(\ref{leadingordergeneral}).
By plugging in (\ref{Ccoefficientsr}), we have that
\begin{eqnarray*}
&&\mathbb{E}\bigl[F_{j_{1}+1,r_{1}}F_{j_{2}+1,r_{2}}\cdots
F_{j_{l}+1,r_{l}}|W(1)=z
\bigr]
\\
&&\qquad=\sum _{ \{
(\mathbf{i}_{1},\mathbf{i}_{2},\ldots,\mathbf{i}_{l})|\llVert\mathbf{i}_{\omega}\rrVert=j_{\omega
}+1,\omega =1,2,\ldots,l \} }\prod_{\omega=1}^{l}C_{\mathbf {i}_{\omega
},r_{\omega}}(x_{0})P_{(\mathbf{i}_{1},\mathbf{i}_{2},\ldots,\mathbf{i}
_{l})} ( z ),
\end{eqnarray*}
where $P_{(\mathbf{i}_{1},\mathbf{i}_{2},\ldots,\mathbf{i}_{l})}(z)$ is
defined in (\ref{ce_product_stratonovichintegrals}).
Formula (\ref{formulageneralOmegak}) follows from the fact that
\begin{eqnarray*}
&&\frac{\partial}{\partial y_{i}}\bigl(P_{(\mathbf{i}_{1},\mathbf
{i}_{2},\ldots,%
\mathbf{i}_{l})} \bigl( \sigma(x_{0})^{-1}D(x_{0})^{-1}y
\bigr) \phi_{\Sigma(x_{0})}(y)\bigr)
\\
&&\qquad= \biggl(\frac{\partial}{\partial y_{i}}P_{(\mathbf{i}_{1},\mathbf
{i}%
_{2},\ldots,\mathbf{i}_{l})} \bigl( \sigma
(x_{0})^{-1}D(x_{0})^{-1}y \bigr)
\\
&&\hspace*{5.2pt}\qquad\quad{}-P_{(\mathbf{i}_{1},\mathbf{i}_{2},\ldots,\mathbf{i}_{l})} \bigl(
\sigma(x_{0})^{-1}D(x_{0})^{-1}y
\bigr) \bigl(\Sigma(x_{0})^{-1}y\bigr)_{i}
\biggr)\phi_{\Sigma(x_{0})}(y)
\end{eqnarray*}
as well as the definition of the differential operators in (\ref
{Doperators}).

\begin{remark}
The above conditioning argument can be justified, when $\partial
^{\mathbf{r}%
}\delta(B(1)-y)$ is regarded as a generalized Wiener functional (random
variable) and the expectation is interpreted in the corresponding
generalized sense as in Watanabe~\cite{WatanabeAnalysisofWiener1987}.
\end{remark}

\subsection{\texorpdfstring{Proof of Theorem \protect\ref{propositionconvergencedensity}}
{Proof of Theorem 2}}

Now, based on Assumption~\ref{aussumptionboundedderivatives}, we introduce
the following uniform upper bounds. For $k\geq1$, let $\mu_{k}$ and $%
\sigma_{k}$ be the uniform upper bounds of the $k$th order derivative
of $%
\mu$ and $\sigma$, respectively, that is,
%
\begin{equation}\label{equationboundforderivatives}
\biggl\llvert\frac{\partial^{(k)}\mu(x;\theta)}{\partial
x^{k}}\biggr\rrvert\leq\mu_{k}\quad\mbox{and}\quad
\biggl\llvert\frac{\partial^{(k)}\sigma
(x;\theta)}{%
\partial x^{k}}\biggr\rrvert\leq\sigma_{k}
\end{equation}
for $(x,\theta)\in\mathbb{R}^{m}\times\Theta$. Also, let $\mu_{0}$
and $\sigma_{0}$ denote the uniform upper bounds of $\llvert\mu
(x_{0};\theta)\rrvert$ and $\llvert\sigma(x_{0};\theta
)\rrvert$ on $(x_{0},\theta)\in K\times\Theta$, respectively,
that is,
%
\begin{equation}\label{equationboundsfordriftandvolcompact}
\bigl\llvert\mu(x_{0};\theta)\bigr\rrvert\leq\mu_{0}
\quad\mbox{and}\quad
\bigl\llvert\sigma(x_{0};\theta)\bigr\rrvert
\leq\sigma_{0}
\end{equation}
for $(x_{0},\theta)\in K\times\Theta$. In order to establish the uniform
convergence in Theorem~\ref{propositionconvergencedensity}, we introduce
the following lemma. When the dependence of parameters is emphasized, we
express $X^{\epsilon}(1)$ as $X^{\epsilon}(1;\theta,x_{0})$ and express
the standardized random variable $Y^{\epsilon}$ defined in %
(\ref{standardizationXY}) as
\[
Y^{\epsilon}(\theta,x_{0})=D(x_{0}) \bigl(
X^{\epsilon}(1;\theta,x_{0})-x_{0} \bigr) /\sqrt{
\Delta}.
\]
In this Appendix, we employ standard notation of Malliavin calculus (see,
e.g., Nualart~\cite{NualartMCbook} and Ikeda and Watanabe \cite
{IkedaWatanabe1989}) and the
theory of Watanabe~\cite{WatanabeAnalysisofWiener1987} and %
Yoshida~\cite{Yoshida1992StatisticsSmallDiffusion,Yoshida1992MLE}. For
the readers'
convenience, a brief survey of some relative theory is provided in the
supplementary material~\cite{autokey73}.

\begin{lemma}
\label{propositionDexpansion} Under Assumption \ref
{aussumptionboundedderivatives}, the pathwise expansion
(\ref{Yexpansion}) holds in the sense of $D^{\infty}$ uniformly in
$(x_{0},\theta)\in K\times\Theta$, that is,
\[
\Biggl\llVert Y^{\epsilon}(\theta,x_{0})-\sum
_{k=0}^{J}\frac
{1}{k!}\,
\frac{\partial^{(k)}Y^{\epsilon}(\theta,x_{0})}{\partial\epsilon
^{k}}%
\bigg\vert_{\epsilon=0}\epsilon^{k}\Biggr
\rrVert_{D_{p}^{s}}=\mathcal{O}%
\bigl(\epsilon^{J+1}
\bigr)
\]
for any $J\in\mathbb{N}$, $p\geq1$ and $s\in\mathbb{N}$.
\end{lemma}

\begin{pf}
See the supplementary material~\cite{autokey73}.
\end{pf}

Because of Assumption~\ref{assumptionpositivedefinite}, Theorem
3.4 in Watanabe~\cite{WatanabeAnalysisofWiener1987} guarantees the
uniform nondegenerate condition, that is,
\[
\limsup_{\epsilon\rightarrow0}\mathbb{E}\bigl[\det\bigl(\Sigma
\bigl(Y^{\epsilon
}(\theta,x_{0})\bigr)\bigr)^{-p}\bigr]<
\infty\qquad\mbox{for any }p\in(0,+\infty).
\]
Let $\Lambda=\mathbb{R}^{m}$ denote a set of indices. For any $%
y=(y_{1},\ldots,y_{m})\in\Lambda$, let us consider a generalized function
defined as $T_{y}(z):=\delta(z-y)$, which is a Schwartz distribution,
that is, $%
T_{y}\in\mathcal{S}^{\prime}(\mathbb{R}^{m})$. Applying Theorem
2.3 in %
Watanabe~\cite{WatanabeAnalysisofWiener1987} and Theorem 2.2 in %
Yoshida~\cite{Yoshida1992MLE}, we obtain that
$T_{y}(Y^{\epsilon}(\theta,x_{0}))$ admits the following asymptotic
expansion: for any arbitrary $J\in\mathbb{N}$,
\[
\delta\bigl(Y^{\epsilon}(\theta,x_{0})-y\bigr):=\sum
_{k=0}^{J}\Phi_{k,(\theta,x_{0})}(y)
\epsilon^{k}+\mathcal{O}\bigl(\epsilon^{J+1}\bigr)\qquad \mbox{in }
D^{-\infty},
\]
uniform in $y\in\Lambda$, $%
x_{0}\in K$ and $\theta\in\Theta$. Here, the correction term $\Phi
_{k,(\theta,x_{0})}(y)$ is given in (\ref{Phikgeneral}). Therefore, we
obtain that
\[
\sup_{y\in\Lambda,x_{0}\in K,\theta\in\Theta}\Biggl\llvert\mathbb{E}%
\Biggl( \delta
\bigl(Y^{\epsilon}(\theta,x_{0})-y\bigr)-\sum
_{k=0}^{J}\Phi_{k,(\theta,x_{0})}(y)
\epsilon^{k} \Biggr) \Biggr\rrvert=\mathcal{O}%
\bigl(
\epsilon^{J+1}\bigr).
\]
%
Hence, by taking into account the transform
(\ref{standardizationXY}), we obtain that
\begin{eqnarray*}
&&
\mathop{\sup_{(x,x_{0},\theta) }}_{\in E\times K\times
\Theta}%
\Biggl\llvert
\mathbb{E}\delta\bigl(X^{\epsilon}(1)-x\bigr)-\frac{\det
D(x_{0})}{\sqrt{%
\Delta^{m}}}\sum
_{k=0}^{J}\Omega_{k} \biggl(
\frac
{D(x_{0}) (
x-x_{0} ) }{\sqrt{\Delta}} \biggr) \epsilon^{k}\Biggr\rrvert\\
&&\qquad=
\mathcal{O}%
\bigl( \Delta^{({J+1-m})/{2}} \bigr),
\end{eqnarray*}
which yields (\ref{densityapproximationerrorboundequation}).

\subsection{\texorpdfstring{Proof of Proposition \protect\ref{theoremconvergenceAMLE}}
{Proof of Proposition 1}}

For $J\geq m$, let
\[
R^{(J)}(\Delta,x|x_{0};\Theta):=\sup_{\theta\in\Theta
}\bigl|p_{X}(
\Delta,x|x_{0};\theta)-p_{X}^{(J)}(
\Delta,x|x_{0};\theta)\bigr|.
\]
By Theorem~\ref{propositionconvergencedensity}, there exists a
constant $%
C>0$ such that $R^{(J)}(\Delta,x|x_{0};\Theta)\leq C\epsilon
^{J+1-m}$ for
any $x\in E$ and $\epsilon$ sufficiently small. Thus, for any positive
integer~$k$, it follows that
\[
\mathbb{E} \bigl[ \bigl|R^{(J)}\bigl(\Delta,X(t+\Delta)|X(t);\Theta
\bigr)\bigr|^{k}|X(t)=x_{0}%
\bigr] \leq C^{k}
\epsilon^{k(J+1-m)}\rightarrow0\qquad\mbox{as }\epsilon\rightarrow0.
\]
By the Chebyshev inequality, $R^{(J)}(\Delta,X(t+\Delta)|X(t);\Theta)$
converges to zero in probability given $X(t)=x_{0}$, that is, for any $%
\varepsilon>0$,
\[
\mathbb{P} \bigl[ \bigl|R^{(J)}\bigl(\Delta,X(t+\Delta)|X(t);\Theta\bigr)\bigr|>
\varepsilon|X(t)=x_{0} \bigr] \rightarrow0\qquad\mbox{as }\epsilon
\rightarrow0.
\]
By conditioning, it follows that
\begin{eqnarray*}
&&\mathbb{P} \bigl[ \bigl|R^{(J)}\bigl(\Delta,X(t+\Delta)|X(t);\Theta
\bigr)\bigr|>\varepsilon%
\bigr]
\\
&&\qquad=\int_{R}\mathbb{P} \bigl[ \bigl|R^{(J)}\bigl(
\Delta,X(t+\Delta)|X(t);\Theta\bigr)\bigr|>\varepsilon|X(t)=x_{0} \bigr]
\mathbb{P}\bigl(X(t)\in dx_{0}\bigr).
\end{eqnarray*}
Because of the fact that
\[
0\leq\mathbb{P} \bigl[ \bigl|R^{(J)}\bigl(\Delta,X(t+\Delta)|X(t);\Theta
\bigr)\bigr|>\varepsilon|X(t)=x_{0} \bigr] \leq1
\]
and $\int_{R}\mathbb{P}(X(t)\in dx_{0})=1$, it follows from the Lebesgue
dominated convergence theorem that
\[
\mathbb{P} \bigl[ \bigl|R^{(J)}\bigl(\Delta,X(t+\Delta)|X(t);\Theta\bigr)\bigr|>
\varepsilon%
\bigr] \rightarrow0\qquad\mbox{as }\epsilon\rightarrow0,
\]
that is,
\[
\mathbb{P} \Bigl[ \sup_{\theta\in\Theta} \bigl\vert p_{X}\bigl(
\Delta,X(t+\Delta)|X(t);\theta\bigr)-p_{X}^{(J)}\bigl(
\Delta,X(t+\Delta)|X(t);\theta\bigr)%
\bigr\vert>\varepsilon\Bigr]
\rightarrow0
\]
as $\epsilon\rightarrow0$. Now, we obtain that
%
\begin{equation}\label{convinprobdensity}
p_{X}^{(J)}\bigl(\Delta,X(t+\Delta)|X(t);\theta
\bigr)-p_{X}\bigl(\Delta,X(t+\Delta)|X(t);\theta\bigr)\stackrel{
\mathbb{P}} {\rightarrow}0
\end{equation}
as $\epsilon\rightarrow0$ uniformly in $\theta\in\Theta$. Following
similar lines of argument as those in the proof of Theorem 2 in %
A{\"{\i}}t-Sahalia~\cite{AitSahalia2002Econometrica} and Theorem 3 in
A{\"{\i}}t-Sahalia~\cite{AitSahalia2008AS}, we arrive at
\[
L_{i}^{(J)}(\theta)\stackrel{\mathbb{P}} {\rightarrow}
L_{i}(\theta)\qquad\mbox{as }\epsilon\rightarrow0\mbox{ uniformly in }
\theta\in\Theta
\]
by the convergence in (\ref{convinprobdensity}) and continuity of
logarithm. Hence, for any arbitrary $n>0$, one obtains the
convergence of log-likelihood $\ell_{n}^{(J)}(\theta
)\stackrel{\mathbb{P}}{\rightarrow} \ell_{n}(\theta)$ uniformly in
$\theta$. Finally, the convergence of $\widehat{\theta}{}^{(J)}_{n}-\widehat{\theta}%
_{n}\stackrel{\mathbb{P}}{\rightarrow}0$ as $\epsilon\rightarrow0$
follows directly
from Assumption~\ref{assumption3timesdifferentiable} and the standard
method employed in %
A{\"{\i}}t-Sahalia~\cite
{AitSahalia2002Econometrica,AitSahalia2008AS}.
\end{appendix}

\section*{Acknowledgements}

I am very grateful to Professor Peter B\"{u}hlmann (Co-Editor), the
Associate Editor and three anonymous referees for the constructive
suggestions. I also thank Professors Yacine A\"{i}t-Sahalia, Mark
Broadie, Song Xi Chen, Ioannis Karatzas, Per Mykland, Nakahiro Yoshida
and Lan Zhang for helpful comments.

\begin{supplement}
\stitle{Maximum-likelihood estimation for diffusion processes via
closed-form density expansions---Supplementary material\\}
\slink[doi]{10.1214/13-AOS1118SUPP} 
\sdatatype{.pdf}
\sfilename{aos1118\_supp.pdf}
\sdescription{This supplementary material contains (1) closed-form
formulas for $%
\Omega_{1}(y),\Omega_{2}(y)$ and $\Omega_{3}(y)$, (2) closed-form
expansion formulas for the examples, (3) detailed plots of errors for
the examples, (4) simulation results for the DMROU and SQR models, (5)
an alternative exhibition of the simulation results, (6) two more
examples for simulation study, (7) a brief survey of the Malliavin
Calculus and Watanabe--Yoshida Theory and (8) proof of Lemma
\ref{propositionDexpansion}.}
\end{supplement}

%

\printaddresses

\end{document}